\documentclass[article]{siamart1116}

\usepackage{amssymb,amsmath,verbatim,bm,color,times,mathtools,multicol}
\usepackage{amsopn}
\usepackage{rotating}

\usepackage{xfrac}
\usepackage[normalem]{ulem}
\usepackage{multirow}
\usepackage[center,small]{caption2}
\usepackage{psfrag}
\usepackage{algorithm,algorithmic}
\usepackage{lipsum}

\usepackage{graphicx,epstopdf}
\usepackage[caption=false]{subfig}

\newcommand{\MATLAB}{\textsc{Matlab}}

\title{a radial basis function based optimization algorithm with \\  regular simplex set geometry in ellipsoidal trust-regions\thanks{Submitted to the editors May 29\textsuperscript{th}, 2018}. This research was supported by the research project EMODO, Flanders Make, the strategic research centre for the manufacturing industry.}

\author{Tom Lefebvre \footnotemark[2]\ \footnotemark[3]
\and Frederik De Belie \footnotemark[2]\ \footnotemark[3]
\and Guillaume Crevecoeur\footnotemark[2] \footnotemark[3]}

\begin{document}
	
	\maketitle
	
	\renewcommand{\thefootnote}{\fnsymbol{footnote}}
	
	\footnotetext[2]{Department of Electrical Energy, Metals, Mechanical Constructions \& Systems. Faculty of Engineering and Architecture. Ghent University. Belgium. Corresponding author: \email{tom.lefebvre@ugent.be}}
	\footnotetext[3]{The authorse are associated to EEDT Decision and Control, Flanders Make, Belgium.}
	
	\renewcommand{\thefootnote}{\arabic{footnote}}
	
	\begin{abstract}
	We present a novel derivative-free interpolation based optimization algorithm. A trust-region method is used where a surrogate model is realised via an interpolation framework. The framework for interpolation is provided by Universal Kriging. A first contribution focusses on the development of an original sampling strategy. A valid model is guaranteed by maintaining a well-poised subset that exhibits the regular simplex geometry approximately. It follows that this strategy improves the scattering of points with respect to the state-of-the-art and, even importantly, assures that the surrogate model exhibits curvature. A second contribution focusses on the generalisation of the spherical trust-region geometry to an ellipsoidal geometry, that to account for local anisotropy of the objective function and to improve the interpolation conditions as seen from the output space. 
	The ensemble method is validated against its direct competitors on a set of multidimensional problems. 
	\end{abstract}
	
	\begin{keywords}
		derivative-free optimization, radial basis functions, trust-region methods, nonlinear optimization
	\end{keywords}
	
	
	\section{Introduction}
	We address the unconstrained minimization of objective function $f:\mathcal{X}\subset\mathbb{R}^n\rightarrow\mathbb{R}$. We assume that $f$ is continuously differentiable, smooth and bounded from below and that $\nabla f$ is Lipschitz continuous on subset $\mathcal{X}$. We restrict the minimization to a bounded convex subset $\mathcal{X}\subset\mathbb{R}^n$.
	\begin{equation*}
	\min_{\boldsymbol{x}\in\mathcal{X}\subset\mathbb{R}^n} ~f(\boldsymbol{x})
	\end{equation*}
	Such a convex subset can be defined through a number of inequality constraints. Yet we assume that none of these constraints will be active at the minimum. The subset $\mathcal{X}$ is thus mainly introduced to search for a local minimum within a range of practical interest. This may also be a means to delimit the search space to the definition domain of the function. To motivate the particular development of the algorithms presented we premise additional assumptions on the computational tractability of the objective function $f$. 
	
	As is the case in many engineering design applications, a single evaluation of the objective is assumed to be associated with the execution of a complex computer simulation or legacy code. Therefore, we assume it to be computationally cumbersome. 
	Secondly, given the high degree of realism that is pursued by modern day simulations and the great level of expert knowledge that is required to fully comprehend the underlying physics, once such simulation codes reach the optimization stage, quite often they are treated as black boxes. This implies that little information about the inner working remains accessible except that a single simulation is known to be a time-consuming endeavour  where sensitivity analyses need to be retrieved by examining input-output behaviour, hence also implying that derivative information is not accessible and becomes intractable to directly approximate.
	
	Altogether we thus assume that the optimization needs to be carried out solely depending on direct evaluations of the objective function itself. The class of Derivative Free Optimization (DFO) methods is tailored exactly to problems that exhibit such challenging computational prerequisites. Among the most popular members of this class are the Nelder-Mead simplex method \cite{nelder1965simplex} and pattern search methods \cite{trosset1997numerical}. This article however focusses on the particular subclass of interpolation based optimization (IBO) methods \cite{wild2008orbit,oeuvray2009boosters,Gutmann2001,powell2006newuoa}.
	
	Given the presumed lack of derivative information, DFO methods are forced to assess the behaviour of the optimization hypersurface indirectly. IBO methods rely on a surrogate model which is computationally simple to evaluate. Such surrogate models can be formed by interpolating a set of sufficiently scattered interpolation points. In our framework we use Universal Kriging as a generalisation of radial basis function (RBF) models \cite{DACE,powell1999recent,sacks1989design}. A brute force approach involves sampling the entire feasible space \cite{crombecq2011novel}, 
	whereas limiting the surrogate modelling space to a compact neighborhood at each iterate of the IBO can reduce the overall computational time \cite{alexandrov1998trust}. Trust-region algorithms have illustrated the benefits of following such a strategy.
	
	
	When the surrogate satisfies Taylor-like error bounds on its function value and gradient for every iteration, convergence to first-order optimality can be guaranteed \cite{conn2009global}. In the context of interpolation, such error bounds can be established depending on the geometry of the used interpolation points \cite{Conn2008}. The geometry of points can be understood as the mutual scattering of points with respect to a certain region. The quality of a geometry can be quantified by a set dependent constant. An arbitrary set of points is well-poised for interpolation within the specified region if the corresponding quality measure turns out to be sufficiently large. 
	
	In order to guarantee a valid interpolation model suited for optimization within a specified trust-region, IBO methods maintain a subset of the available interpolation set for which the geometrical constant is maximized. Consequently, the well-poisedness of the subset geometry is assured and so is the model validity. In every iteration this constant is maximized by removing points from the interpolation set and generating new ones according to the current size and location of the trust-region; that whilst limiting the overall accumulation of objective evaluations. This delicate mechanism is referred to as the interpolation set management.
	
	In \cite{oeuvray2009boosters,conn1997convergence,powell2001lagrange} such a quality measure constant is described through a basis of polynomials, other authors have described this constant in function of the QR pivots of the Vandermonde matrix that depends on the interpolation points \cite{wild2008orbit,Conn2008}. All the mentioned methods consider a subset of $n+1$ interpolation points to guarantee a valid model, which they complement with interpolation points from previous iterations if the expanded geometry does not become ill-conditioned \cite{wild2008orbit}. In any of these cases, the interpolation set management implies orthogonality between the $n$ interpolation points with respect to the current trust-region centre. The centre then serves as the $(n+1)^\text{th}$ interpolation point.
	
	
	
	
	In light of this framework, the contribution of this article is twofold. 
	Firstly, in our algorithm we rely on a subset of $n+2$, instead of $n+1$ interpolation points as in \cite{wild2008orbit,oeuvray2009boosters}. This way, we aim at achieving an improved scattering which is heuristically justified for low dimensional problems. Considering the RBF interpolation framework used, this also means that at every iteration the surrogate will exhibit curvature which is of vital essence for the well-functioning of our second contribution. This first, seemingly simple, modification has a large impact on the interpolation set management. An original sampling strategy is devised that is inspired by the volumetric extremity properties of the regular $n$-simplex. Details on the exact motivation and consequences are provided in the accompanying sections.
	
	A second contribution focusses on a generalisation of the conventional spherical trust-region framework to an ellipsoidal trust-region geometry. The motivation to use an ellipsoidal geometry is again twofold. An ellipsoidal trust-region will adapt to the geometry of the objective function and will therefore allow larger steps towards more promising regions compared to using a spherical geometry. We will inspect the convergence properties by means of extensive numerical experiments. Secondly our interpolation set management is input space oriented. The introduction of an ellipsoidal trust-region, which shape reflects the objective's local output behaviour, allows to indirectly account for the output behaviour when executing the interpolation set management and will benefit the interpolation based metamodelling.

	\subsection{Outline}
	To keep this article self-contained, we shortly resume the general trust-region framework for surrogate management in \cref{sec:trustframe}. In \cref{sec:RBF}, we elaborate on the radial basis function interpolation framework and provide details on the implementation. Section \ref{sec:fullylin} motivates our choice for using $n+2$ interpolation points and introduces our interpolation set management. In \cref{sec:ellips}, we discuss the transition towards ellipsoidal trust regions. In \cref{sec:numericalresults}, we benchmark the proposed methods using extensive numerical experiments and in \cref{sec:conclusions} conclusions are drawn.

	\section{Trust-region methods}
	\label{sec:trustframe}
	In this section we shortly review the trust-region methodology for DFO in a generic manner. Furthermore, we define fully linear surrogate models that generalise the characteristics of the surrogate model used by the trust-region methodology from derivative-based models to more general model types. As made apparent in the introduction, we are primarily interested in the use of interpolation models. Finally, we reach out towards the inclusion of the trust-region framework in the context of IBO.  
	
	\subsection{Trust-region framework}
	The trust-region framework is a widely-used technique to attain global convergence of derivative-based optimization algorithms. 
	At every iteration $k$, a derivative-based surrogate model, $m_k$, is considered that approximates the actual objective within a neighbourhood of iterate $\boldsymbol{x}_k$. Traditionally, this neighbourhood, i.e. the trust-region, is modelled as the parameterized set,
	\begin{equation}
	\mathcal{B}(\boldsymbol{x}_k,\Delta_k) = \left\lbrace\boldsymbol{x}\in\mathbb{R}^n:\|\boldsymbol{x}-\boldsymbol{x}_k\|\leq\Delta_k\right\rbrace
	\end{equation}
	Here $\boldsymbol{x}_k$ is the current iterate and $\Delta_k$ is referred to as the trust-region radius. The choice of the norm is arbitrary, commonly the standard Euclidean distance is used. In that case the trust-region becomes an $n$-dimensional hypersphere.
	
	A trial step, $\boldsymbol{s}_k$, is obtained by minimizing surrogate, $m_k$, in the trust-region, $\mathcal{B}(\boldsymbol{x}_k,\Delta_k)$. This minimization problem is referred to as the trust-region subproblem. The resulting trial step, $\boldsymbol{s}_k$, gives rise to the trivial candidate iterate, $\boldsymbol{x}_{k+1} = \boldsymbol{x}_k + \boldsymbol{s}_k$. 
	\begin{equation}
	\label{eq:TR:subproblem}
	\boldsymbol{s}_k = \arg \min_{\boldsymbol{s}\in\mathcal{B}(\boldsymbol{0},\Delta_k)} m_k(\boldsymbol{x}_k+\boldsymbol{s})
	\end{equation}

	In between iterations, the size and position of the trust-region is updated such that the generated sequence of iterates converges. The configuration of the trust-region, which is entirely determined by the couple $(\boldsymbol{x}_k,\Delta_k)$, is updated according to a fundamental scheme based on the ratio of actual to predicted decrease, $\rho_k$. The ratio constitutes a metric for the predictive quality of the surrogate with respect to the current trust-radius \cite{alexandrov1998trust}.
	\begin{equation}
	\label{eq:TR:ratio}
	\rho_k  = \frac{f(\boldsymbol{x}_k)-f(\boldsymbol{x}_k+\boldsymbol{s}_k)}{m_k(\boldsymbol{x}_k)-m_k(\boldsymbol{x}_{k}+\boldsymbol{s}_k)} 
	\end{equation}
	
	Based on the value of $\rho_k$, the trust region is updated according to the following general rules. If the prediction of the surrogate is good, the radius is increased; in case of bad prediction, the radius is reduced. For other cases, the current radius is simply maintained.
	
	To regulate the expansion and contraction of the trust-region, one takes values $0<\Delta_m\leq\Delta_M,0\leq\eta_1<\eta_2<1$ and $0<\gamma_1<1<\gamma_2$. Here $\Delta_m$ and $\Delta_M$ represent a minimal and maximal allowed trust-region radius, respectively. The radius is then updated according to
	\begin{equation}
	\label{eq:TR:radiusupdate}
	\Delta_{k+1} = \left\lbrace\begin{aligned}
	& \max\left\{\gamma_1\Delta_k,\Delta_m \right\},& &\rho_k\leq\eta_1\\
	& \Delta_k,& &\eta_1<\rho_k\leq\eta_2 \\
	& \min\left\{\gamma_2\Delta_k,\Delta_M\right\},& & \eta_2 < \rho_k
	\end{aligned}\right.
	\end{equation}
	
	To decide whether to accept or reject a trial step, $\boldsymbol{s}_k$, we adopt the rules from \cite{alexandrov1998trust}. If the candidate iterate, $\boldsymbol{x}_k+\boldsymbol{s}_k$, is not accepted then the trust-radius shrinks given that $\rho_k < 0$.
	\begin{equation}
	\label{eq:TR:iterateupdate}
	\boldsymbol{x}_{k+1} = \left\lbrace\begin{aligned}
	&\boldsymbol{x}_k + \boldsymbol{s}_k ,& & f(\boldsymbol{x}_k+\boldsymbol{s}_k)<f(\boldsymbol{x}_k) \\
	&\boldsymbol{x}_k ,& & \text{otherwise}
	\end{aligned}\right.
	\end{equation}
	
	
	\begin{algorithm}[t]
		\caption{Generic trust-region algorithm}
		\label{alg:TR:general}
		\begin{algorithmic}[1]
			\STATE{choose $\boldsymbol{x}_0\in\mathcal{X}\subset\mathbb{R}^n$, $0<\Delta_m\leq\Delta_0\leq\Delta_M$, $0\leq\eta_1<\eta_2<1$ and $0<\gamma_1<1<\gamma_2$}
			\FOR{$k = 0,1,\dots$ until convergence}
			\STATE{obtain model $m_k$}
			\STATE{solve subproblem \cref{eq:TR:subproblem} to obtain $\boldsymbol{s}_k$}
			\STATE{evaluate $f(\boldsymbol{x}_k + \boldsymbol{s}_k)$ and compute $\rho_k$ through \cref{eq:TR:ratio}}
			\STATE{update $\boldsymbol{x}_{k+1}$ according to \cref{eq:TR:iterateupdate} and $\Delta_{k+1}$ according to \cref{eq:TR:radiusupdate} }
			\ENDFOR
			\RETURN $\boldsymbol{x}_k$ 
		\end{algorithmic} 
	\end{algorithm}
	
	The generic trust-region algorithm is presented in \cref{alg:TR:general}. 

	\subsection{Trust-region method using fully linear models} In the derivative-based version of the trust-region method, the surrogate, $m_k$, is provided by a quadratic model of the form
	\begin{equation}
	\label{eq:tr:trmodel}
	m_k(\boldsymbol{x}) = f(\boldsymbol{x}_k) + \boldsymbol{g}_k^\top (\boldsymbol{x}-\boldsymbol{x}_k) + \frac{1}{2} (\boldsymbol{x}-\boldsymbol{x}_k)^\top {\boldsymbol{\mathrm{H}}}_k  (\boldsymbol{x}-\boldsymbol{x}_k)
	\end{equation}
	Here, $\boldsymbol{g}_k$ and ${\boldsymbol{\mathrm{H}}}_k$ are evaluations of the first- and second-order derivatives of the objective function $f$ in the iterate $\boldsymbol{x}_k$, respectively. According to Taylor's theorem, a region exists within the vicinity of $\boldsymbol{x}_k$ where the assumption that model $m_k$ is a reliable approximation of $f$, is true. A region of reliable approximation can be defined as a region where the surrogate modelling approximation error, reflected by $\rho_k$, is small enough such that descend in the surrogate model, $m_k$, restricted to this region corresponds to descend of the objective fuction $f$. It is apparent that the strength of the trust-region method thus depends on the capacity of the updating rules to maintain a region where the surrogate modelling approximation error is sufficiently bounded, so that $m_k$ is suited for optimization purposes \cite{conn2000trust,dennis1997convergence}.
	
	From that perspective the generic trust-region method can be decomposed into two crucial mechanisms that interact to generate the sequence of iterates. A procedure must be provided that can generate a quadratic approximation model $m_k$ of the goal function $f$ for any given iterate $\boldsymbol{x}_k$. The updating rules will then maintain a series of regions where the corresponding and consecutive models, $m_k$ are suited to generate a descending sequence of iterates as the solution of subproblem \cref{eq:TR:subproblem}. The latter mechanism is in that aspect generic and can be employed to maintain a series of regions given any type of model, also derivative-free.
	
	In this work we are interested in the class of so called fully linear models. The definition of a fully linear model establishes Taylor-like bounds on arbitrary models and therefore generalises the notion of a reliable approximation. As a consequence, the trust-region framework can be generalised to include such fully linear models. 
	
	We adopt the definition from \cite{conn2009global,Conn2008} with $\hat{\mathcal{X}} = \{\boldsymbol{x}\in\mathcal{X} : \mathcal{B}(\boldsymbol{x},\Delta_M) \subseteq \mathcal{X}\}$.
	\begin{definition}
		\label{def:fullylinear}
		For given $\kappa_f>0$, $\kappa_g>0$, and given $\boldsymbol{x}_k\in\hat{\mathcal{X}}$, $\Delta_k\in(0,\Delta_M]$ defining $\mathcal{B}_k = \mathcal{B}(\boldsymbol{x}_k,\Delta_k)$, surrogate model $m_k$ is said to be fully linear on $\mathcal{B}_k$, if for all $\boldsymbol{x}\in \mathcal{B}_k$:
		\begin{subequations}
		\begin{eqnarray}
			|f(\boldsymbol{x}) - m_k(\boldsymbol{x})| &\leq \kappa_f \Delta_k^2 \label{eq:fullylinear:a}\\
			\|\nabla f(\boldsymbol{x}) - \nabla m_k(\boldsymbol{x})\| &\leq \kappa_g \Delta_k \label{eq:fullylinear:b}
		\end{eqnarray}
		\end{subequations}
	\end{definition}
	
	
	Let us assume two procedures that can either generate a fully linear surrogate model $m_k$ on $\mathcal{B}(\boldsymbol{x}_k,\Delta_k)$ on the one hand, or that can verify whether or not the surrogate model, $m_k$, is fully linear on $\mathcal{B}(\boldsymbol{x}_k,\Delta_k)$ on the other, and that for any given $\boldsymbol{x}_k\in\hat{\mathcal{X}}$ and $\Delta_k \in (0,\Delta_M]$. If two such procedures are available, \cref{alg:TR:general} can be modified to construct a trust-region algorithm with fully linear models that will converge to a first-order critical point of the function $f$ \cite{conn2009global}. \cref{alg:TR:fullylinear} presents such a trust-region algorithm using fully linear models.
	
	The main modification is the introduction of the termination criteria, $\| \nabla m_k(\boldsymbol{x}_k) \| \leq \frac{\epsilon_g}{2}$. This condition is a reformulation of the first-order optimality condition on the objective function, $\| \nabla f(\boldsymbol{x}_k) \| \leq {\epsilon_g}$. Assuming that $m_k$ is fully linear on $\mathcal{B}(\boldsymbol{x}_k,\Delta_k^g)$, we have that
	\begin{equation}
	\begin{aligned}
	\|\nabla f(\boldsymbol{x}_k)\| &\leq \|\nabla m_k(\boldsymbol{x}_k)\| + \|\nabla f(\boldsymbol{x}_k)-\nabla m_k(\boldsymbol{x}_k)\| \\ &\leq \frac{\epsilon_g}{2} + \kappa_g \Delta^g_k
	\end{aligned}
	\label{eq:tr:firstorderopt}
	\end{equation}
	
	\cref{alg:TR:fullylinear} assures that upon termination, model $m_k$ will be fully linear on $\mathcal{B}(\boldsymbol{x}_k,\frac{\epsilon_g}{2\kappa_g})$. Substituting $\frac{\epsilon_g}{2\kappa_g}$ for $\Delta^g_k$ in \cref{eq:tr:firstorderopt} reveals that this is equivalent with satisfying the first-order optimality criterion, $\|\nabla f(\boldsymbol{x}_k)\|\leq \epsilon_g$. We remark that it is thus fundamental to provide procedures that can guarantee the tightest possible bound on $\kappa_f$ and $\kappa_g$.	In \cref{sec:fullylin}, such procedures are developed considering that model $m_k$ is an interpolation model.
	
	
	\begin{algorithm}[t]
		\caption{Trust-region algorithm with fully linear models}
		\label{alg:TR:fullylinear}
		\begin{algorithmic}[1]
			\STATE{choose $\boldsymbol{x}_0\in\mathcal{X}\subset\mathbb{R}^n$, $0<\Delta_m\leq\Delta_0\leq\Delta_M$, $0\leq\eta_1<\eta_2<1$ and $0<\gamma_1<1<\gamma_2$}
			\FOR{$k = 0,1,\dots$ }
			\STATE{construct fully linear model $m_k$ valid in  $\mathcal{B}(\boldsymbol{x}_k,\Delta_k)$}
			\IF{$\|\nabla m_k(\boldsymbol{x}_k)\|\leq \frac{\epsilon_g}{2}$}
			\IF{$m_k$ is fully linear in $\mathcal{B}(\boldsymbol{x}_k,\frac{\epsilon_g}{2\kappa_g})$}
			\STATE{stop}
			\ELSE
			\STATE{set $\boldsymbol{x}_{k+1}=\boldsymbol{x}_k$ and $\Delta_{k+1} = \frac{\epsilon_g}{2\kappa_g}$}
			\ENDIF
			\ELSE
			\STATE{solve subproblem \cref{eq:TR:subproblem} to obtain $\boldsymbol{s}_k$}
			\STATE{evaluate $f(\boldsymbol{x}_k + \boldsymbol{s}_k)$ and compute $\rho_k$ through \cref{eq:TR:ratio}}
			\STATE{update $\boldsymbol{x}_{k+1}$ according to \cref{eq:TR:iterateupdate} and $\Delta_{k+1}$ according to \cref{eq:TR:radiusupdate}}
			\ENDIF
			\ENDFOR
			\RETURN $\boldsymbol{x}_k$ 
		\end{algorithmic} 
	\end{algorithm}
	
	\subsection{Interpolation models}
	\label{sec:tr:interpol}
	When first-order derivative information, $\boldsymbol{g}_k$, is not directly available, it is possible to construct a surrogate model based on strategically generated data, i.e. direct evaluations of the objective function. A model $m_k$ denotes an interpolation model if it satisfies the interpolation conditions in \cref{eq:interpolationcondition}, for a given input set, $\mathcal{Z} = \{\boldsymbol{z}_1,\dots,\boldsymbol{z}_q\}$ with interpolation points $\boldsymbol{z}_i\in\mathbb{R}^n$, 
	and corresponding output set, $\mathcal{F} = \{f(\boldsymbol{z}_1),\dots,f(\boldsymbol{z}_q)\}$.
	\begin{equation}
	\label{eq:interpolationcondition}
	m_k(\boldsymbol{z}_i) = f(\boldsymbol{z}_i), ~ i = 1,\dots,q
	\end{equation}
	
	
	A quadratic model, as in \cref{eq:tr:trmodel}, is capable of modelling the curvature of the underlying function $f$ and allows therefore to proof global convergence to second-order critical points \cite{conn2009global}. Correspondingly, several DFO methods have been developed that operate with quadratic interpolation models \cite{powell2006newuoa,powell2002uobyqa,conn1997recent}. To fully determine a quadratic model, the number of interpolation points, $q$, must either equal or exceed ${q} = \frac{1}{2}(n+1)(n+2)$ requiring thus at least that number of objective function evaluations. This requirement contradicts with our intention to limit the number of objective function evaluations. 
	
	A linear model can be constructed with as little as $n+1$ interpolation points. Algorithms that operate with minimal interpolation sets of ${q} = n+1$ are described by \cite{wild2008orbit,oeuvray2009boosters,powell1994direct}. The drawback in this case is, however, that a linear model is not capable of representing any curvature present in the objective. 
	In this work we impose the additional requirement that the minimal number of interpolation points must equal ${q} = n+2$, because it is then guaranteed that the resulting model can represent curvature given a proper interpolation framework (see \cref{sec:RBF}). We will refer to such models as \textit{better-than-linear} models, noting that these models classify as fully linear models but are constructed with at least one more interpolation point than the minimum for exact linear models $q=n+1$. 
	These better-than-linear models are further detailed and engaged in \cref{sec:fullylin} and \cref{sec:ellips}.
	
	\subsection{Interpolation based optimization} Anticipating that a procedure is available that can generate an interpolation model, $m_k$, that classifies as a fully linear model, an interpolation framework can be embedded in \cref{alg:TR:fullylinear}. 
	
	Section \ref{sec:RBF} describes a multivariate interpolation framework that can model curvature and is flexible in the sense that it does not impose strict requirements on the number of interpolation points. Within the context of this framework, it will be shown that an interpolation model can be classified as a fully linear model if the set, $\mathcal{Z}$, satisfies conditions on the geometry of the points with respect to the trust-region \cite{Conn2008}. The geometry can be understood as the scattering of interpolation points inside a specified region. In \cref{sec:fullylin} the geometry of an interpolation set is quantified through a single set dependent constant that makes it possible to classify an interpolation model as a fully-linear model.
	
	The procedures of verifying and generating fully linear interpolation models thus boil down to verifying whether the geometric conditions are satisfed with respect to having fully-linear model for a given interpolation set and to complement an existing or possibly empty set so that the condition is satisfied by construction. In the context of IBO, these two mechanisms are provided by a single procedure that is referred to as the interpolation set management or the sampling strategy. A sampling strategy for minimal interpolation sets of ${q}=n+2$ is proposed in \cref{sec:fullylin}. 
	
	
	
	\section{Radial Basis Functions}
	\label{sec:RBF}
	In IBO, the previously described \cref{alg:TR:fullylinear} embeds interpolation models \cref{eq:interpolationcondition}. This requires a flexible, multivariate interpolation modelling framework that can capture curvature and does not impose strict conditions on the number of interpolation points since we will recycle interpolation points (and their evaluations) from previous iterations. The Radial Basis Function (RBF) interpolation framework has readily shown its utility within IBO \cite{wild2008orbit,oeuvray2009boosters,Gutmann2001}. The traditional RBF framework decides arbitrarily on the shape of the modelling RBFs by fixing a shape parameter a priori. Alternatively one can determine the shape parameter based on the interpolation data. In that case the RBF framework coincides with Universal Kriging \cite{DACE} where the shape parameter is referred to as a hyperparameter. The modelling assumptions of the framework and the determination of the model- and hyperparameters are discussed next.
	
	\subsection{Modelling framework}
	Consider a given set of interpolation points $\mathcal{Z}$ and corresponding objective function evaluations $\mathcal{F}$. The RBF multivariate interpolation framework then presumes a function model of the form 
	\begin{equation}
	m(\boldsymbol{x}) = \sum_{\boldsymbol{z}_i\in\mathcal{Z}} \alpha_i\phi\left(\left\|\boldsymbol{x}-\boldsymbol{z}_i\right\|\right) + p(\boldsymbol{x})
	\end{equation}
	where $\phi:\mathbb{R}_+\rightarrow\mathbb{R}$ is an element of the class of RBFs and $p(\boldsymbol{x})$ is an element of an $n$-variate polynomial space $\mathcal{P}^n_d$ of at most degree $d$ referred to as the polynomial tail. Let $\hat{\mathcal{P}}^n_d = \left\{\pi_1(\boldsymbol{x}),\dots,\pi_{\hat{d}}(\boldsymbol{x})\right\}$ be a $n$-variate basis for $\mathcal{P}^n_d$ with $\hat{d}$ the size of this $n$-variate basis, then the polynomial tail can be represented by a linear combination of this basis, that is $p(\boldsymbol{x}) = \sum_{\pi_j \in\hat{\mathcal{P}}^n_d} \beta_j\pi_j(\boldsymbol{x})$ with $\beta_j$ the expansion coefficients. Some popular choices for the RBF are listed in \cref{tab:rbfs}. Notice that the function $\phi(\|\boldsymbol{x}\|)$ has isosurfaces coinciding with hyperspheres in $\mathbb{R}^n$, i.e. the function radiates out isotropically.
	
	Introducing the vectors, $\boldsymbol{\phi}(\boldsymbol{x})$ and $\boldsymbol{\alpha}$, $\in\mathbb{R}^q$ and, $\boldsymbol{\pi}(\boldsymbol{x})$ and $\boldsymbol{\beta}$, $\in\mathbb{R}^{\hat{d}}$ which entries are given by $\phi(\|\boldsymbol{x}-\boldsymbol{z}_i\|),\alpha_i,\pi_j(\boldsymbol{x})$ and $\beta_j$ respectively, the RBF model can be reformulated as
	\begin{equation}
	\label{eq:rbfmodel}
	m(\boldsymbol{x}) =\sum_{\boldsymbol{z}_i\in\mathcal{Z}} \alpha_i\phi\left(\left\|\boldsymbol{x}-\boldsymbol{z}_i\right\|\right)+\sum_{\pi_j \in\hat{\mathcal{P}}^n_d} \beta_j\pi_j(\boldsymbol{x}) = \boldsymbol{\phi}(\boldsymbol{x})^\top\cdot\boldsymbol{\alpha} + \boldsymbol{\pi}(\boldsymbol{x})^\top\cdot\boldsymbol{\beta}
	\end{equation}
	
	\begin{table}[t!]
		\caption{Radial basis functions, $\gamma>0$, where $d_0$ is the minimal polynomial degree to have CPD.}
		\label{tab:rbfs}
		\centering
		\begin{tabular}{|c|c|c|} \hline
			type & $\phi(r)$ & $d_0$\\ \hline\hline
			{cubic} & $r^3$ & 1\\
			{thin plate spline} & $r^2\log(r)$ & 1\\
			{linear} & $\max\{0,1-\gamma|r|\}$ & 0\\
			{multiquadratic} & $\sqrt{r^2+\gamma^2}$ & 0\\
			{gaussian} & $\exp\left(-\gamma r^\beta\right)$ & -1\\ \hline
		\end{tabular}
	\end{table}
	
	\subsection{Determination of model parameters}\label{sec:sub:modelparameters}	The RBF function model (\ref{eq:rbfmodel}) is uniquely defined by the parameters $\boldsymbol{\alpha}$ and $\boldsymbol{\beta}$. These parameters can be determined as follows.
	
	Let the Vandermonde matrices ${\boldsymbol{\Phi}}\in\mathbb{R}^{q\times q}$ and ${\boldsymbol{\Pi}}\in\mathbb{R}^{q\times \hat{d}}$ be defined as $\Phi_{ik} = \phi(\|\boldsymbol{z}_i-\boldsymbol{z}_k\|)$ and $\Pi_{ij} = \pi_j(\boldsymbol{z}_i)$. The set of $q$ interpolation conditions in \cref{eq:interpolationcondition} can then be written as 
	\begin{equation}
	\label{eq:interpol-linsys}
	\begin{bmatrix}
	{\boldsymbol{\Phi}} & {\boldsymbol{\Pi}}
	\end{bmatrix}\cdot\begin{bmatrix}
	\boldsymbol{\alpha} \\ \boldsymbol{\beta}
	\end{bmatrix} = \boldsymbol{f}
	\end{equation}
	Here the entries of $\boldsymbol{f}$ are defined as the function evaluations contained in $\mathcal{F}$.
	
	These conditions do not render unique $\boldsymbol{\alpha}$ and $\boldsymbol{\beta}$. Therefore $\hat{d}$ additional conditions are required to complement \cref{eq:interpol-linsys} and so that the assembled system matrix is nonsingular. To this end we require of $m(\boldsymbol{x})$ that it is the \textit{smoothest} function of the form \cref{eq:rbfmodel} that complies to the interpolation condition \cref{eq:interpolationcondition}. For RBFs, function smoothness can be defined in relation to the concept of \textit{conditional positive definiteness} (CPD) \cite{powell1999recent}. RBFs have the specific property that to any RBF, $\phi$, one can associate a value $d_0$ 
	such that if the maximum degree $d$ of the polynomial tail satisfies $d\geq d_0$, the matrix ${\boldsymbol{\Phi}}$ will be positive definite with respect to all vectors $\boldsymbol{\alpha} \in \mathrm{null}\left({\boldsymbol{\Pi}}\right)$. Equivalently, we have that if $d\geq d_0$ and ${\boldsymbol{\Pi}}^\top \boldsymbol{\alpha} = 0$ it holds that $(-1)^{d_0+1}\boldsymbol{\alpha}^\top{\boldsymbol{\Phi}}\boldsymbol{\alpha} > 0$. In conclusion, when CPD is satified, i.e. $\boldsymbol{\alpha}^\top{\boldsymbol{\Phi}}\boldsymbol{\alpha}>0$, it can be shown that the corresponding interpolator is the smoothest function of the form \cref{eq:rbfmodel} satisfying \cref{eq:interpolationcondition}. For an elaborate proof and additional details on CPD, we refer to \cite{powell1999recent}.
	
	Formally, the matrix representation of the interpolation condition can be complemented as such by the function smoothness condition so to obtain the linear system of equations
	\begin{equation}
	\label{eq:rbf-linsys}
	\begin{bmatrix}
	{\boldsymbol{\Phi}} & {\boldsymbol{\Pi}} \\
	{\boldsymbol{\Pi}}^\top & {\boldsymbol{0}}
	\end{bmatrix}\cdot\begin{bmatrix}
	\boldsymbol{\alpha} \\
	\boldsymbol{\beta}
	\end{bmatrix} = \begin{bmatrix}
	\boldsymbol{f}\\
	\boldsymbol{0}
	\end{bmatrix}
	\end{equation}
	
	\newpage
	The model parameters are determined uniquely if the system matrix is nonsingular. It is easily verified that if $d\geq d_0$, it is then sufficient that ${\boldsymbol{\Pi}}^\top$ has full column rank. If so, the interpolation set is said to be well-poised for interpolation. Such corresponds with geometric conditions on the points in $\mathcal{Z}$ as will be discussed in \cref{sec:fullylin}. 
	
	Based on the previous, the explicit solution of the model parameters is
	\begin{subequations}
		\label{eq:explicitRBFsol}
		\begin{eqnarray}
		\boldsymbol{\beta} &= &\left({\boldsymbol{\Pi}}^\top{\boldsymbol{\Phi}}^{-1}{\boldsymbol{\Pi}}\right)^{-1}{\boldsymbol{\Pi}}^\top{\boldsymbol{\Phi}}^{-1}\boldsymbol{f} \\
		\boldsymbol{\alpha} &= &{\boldsymbol{\Phi}}^{-1}(\boldsymbol{f}-{\boldsymbol{\Pi}}\boldsymbol{\beta}) 
		\end{eqnarray}
	\end{subequations} 
	
	This result may also be obtained through the Universal Kriging (UK) framework. Here, a linear predictor of the form $m(\boldsymbol{x}) = \boldsymbol{\eta}(\boldsymbol{x})^\top \boldsymbol{f}$, is considered. The objective function is modelled as the realization of a polynomial tail and a random process, $f(\boldsymbol{x}) = p(\boldsymbol{x}) + w(\boldsymbol{x})$. The process' stochasticity is quantified by a spatial correlation function that models the output correlation as a function of the spatial distance between input values. The correlation is thus expressed in terms of a RBF as $\mathrm{corr}(w(\boldsymbol{z}_i),w(\boldsymbol{z}_j)) = \phi(\|\boldsymbol{z}_i-\boldsymbol{z}_j\|)$. It can be shown, e.g. \cite{DACE,sacks1989design}, that the optimal linear unbiased predictor, $\boldsymbol{\eta}$, is then determined by the following system of equations where $\boldsymbol{\xi}$ can be interpreted as a Lagrangian multiplier. 
	
	\begin{equation}
	\label{eq:kriging-linsys}
	\begin{bmatrix}
	{\boldsymbol{\Phi}} & {\boldsymbol{\Pi}} \\
	{\boldsymbol{\Pi}}^\top & {\boldsymbol{0}}
	\end{bmatrix}\cdot\begin{bmatrix}
	\boldsymbol{\eta} \\
	\boldsymbol{\xi}
	\end{bmatrix} = \begin{bmatrix}
	\boldsymbol{\phi}\\
	\boldsymbol{\pi}
	\end{bmatrix}
	\end{equation}
	
	For additional details we refer to \cref{appendix:unikrig}.
	
	\subsection{Determination of model hyperparameters}
	\label{sec:sub:hyper}
	The value of the shape parameter $\gamma$ (see \cref{tab:rbfs}) is arbitrary and affects the receding radial effect of the RBF contribution in \cref{eq:rbfmodel}. Heuristics may be provided, however they usually lack proper theoretical foundation. A suitable value is mostly determined by trial-and-error. Furthermore we emphasize that each component of the interpolation points $\boldsymbol{z}_i$ contributes equally to the distance metric, $\|\cdot\|$, i.e. the model presumes isotropic dependence of the objective function on all of its variables. Such isotropic behaviour can be achieved by normalizing the optimization space by performing a linear transformation of the input space. Nonetheless, such measure cannot account for any local deviations on the global trend. Since we are to construct local interpolation models which validity is not to stretch beyond the limits of the trust-region boundary, we could redefine this linear transformation locally. 
	
	To account for local anisotropy in the principle directions of the input space one may define the distance metric as \cref{eq:distance}. The RBFs in \cref{eq:rbfmodel} are then redefined as $\phi(\|\boldsymbol{x}-\boldsymbol{z}_i\|_{\boldsymbol{\gamma}})$. It is clear that in this case the univariate hyperparameter $\gamma$ has become redundant. In practical implementations of the UK framework, the determination of hyperparameters $\gamma_i$ is automated through a maximum log likelihood estimate based on the observations, $\mathcal{F}$. 
	
	\begin{align}
	\label{eq:distance}
	\|\boldsymbol{x}\|^2_{\boldsymbol{\gamma}} =  \sum_{i=1}^n \gamma_i x_i^2
	\end{align}
	
	We refer to \cref{appendix:unikrig} and \cite{toal2009proper} for elaborate detail.
	
	\newpage
	\section{Interpolation set management}
	\label{sec:fullylin}
	In \cref{sec:trustframe}, \cref{alg:TR:fullylinear} described a trust-region method with fully linear models. It was demonstrated how a sequence of fully linear models $m_k$ can be utilised to converge to a first-order critical point of the objective function $f$. When considering interpolation models satisfying \cref{eq:interpolationcondition}, the question whether $m_k$ is a fully linear model boils down to a proper management of the interpolation set $\mathcal{Z}$ used to construct $m_k$. This so-called sampling strategy should maintain a proper geometry of the interpolation set.
	
	In \cite{Conn2008} the geometry of an arbitrary interpolation set, $\mathcal{Z}$, was quantified by a single set dependent constant directly related to the conditions in \cref{eq:fullylinear:a} to define a fully linear model. The lemma in \cite{Conn2008} however, only applies to interpolation models for which the set size equals ${q} = n+1$ exactly. On this basis the authors in \cite{wild2008orbit} devised an IBO algorithm employing RBF interpolation models given at least $n+1$ interpolation points to construct a fully linear model. By applying a QR pivoting strategy directly related to the measure described by \cite{Conn2008}, their sampling strategy assures that a subset of $n+1$ interpolation points (amongst the entire available data set) exhibits the required geometry to guarantee a fully linear interpolation model. Possibly additional samples need to be generated to have verification of the full linearity of the interpolation model. We henceforth refer to these points to verify that $m_k$ is fully linear, as the affine subset. Additional points from the entire data set are subsequently recycled to enhance the modelling capabilities of the RBF model as long as the well-poisedness of \cref{eq:rbf-linsys} is not corrupted.

	Having an interpolation model $m_k$ representing curvature, requires an affine subset of $q\geq n+2$, cfr. \cref{sec:sub:modelparameters}, that is adequately embedded in an interpolation based trust-region algorithm. In other words, we aim to ensure convergence of \cref{alg:TR:fullylinear} to second-order critical points whilst only adding a single point to the affine subset with respect to \cite{wild2008orbit,oeuvray2009boosters}. Since an estimate of the curvature is now always available, the remark on the RBF hyperparameters in \cref{sec:sub:hyper} can be generalised to arbitrary directions. The latter is investigated in \cref{sec:ellips}.
	
	
	In \ref{sec:sub:better-than-linear}, we introduce a generalised version of the theorem in \cite{Conn2008} for interpolation sets with size exceding $n+1$. On that basis we identify the optimal affine subset geometry that corresponds with $n+2$ interpolation points in \ref{sec:sub:optimal-geometry}. Finally in \ref{sec:sub:regular-simplex}, we propose a procedure capable to maintain an interpolation set with an affine subset of $n+2$ points so that interpolation model, $m_k$, is fully linear and exhibits curvature.

	\subsection{Better-than-linear models}
	\label{sec:sub:better-than-linear}
	Uniqueness of the RBF interpolation model parameters implied well-poisedness, i.e. full column rank, of the matrix $\boldsymbol{\Pi}^\top$. The size of the interpolation set must thus equal or exceed the size of the polynomial basis. Given that we adress a minimal interpolation set size $n+2$ and thus $\hat{d}\leq n+2$, the polynomial tail will be linear or equivalently $p(\boldsymbol{x})\in \mathcal{P}_1^n$. A convenient choice for the basis $\hat{\mathcal{P}}^n_{1}$ is then clearly $\boldsymbol{\pi}(\boldsymbol{x})^\top = \begin{bmatrix}\boldsymbol{x}^\top & 1\end{bmatrix}$. Since we will embed the interpolation framework in a trust-region framework, the current iterate $\boldsymbol{x}_k$ is always available for interpolation since $f(\boldsymbol{x}_k)$ is evaluated for calculating the predictive quality of the interpolation model \cref{eq:TR:ratio}. Henceforth, we will consider a trust-region centered at the origin, $\boldsymbol{0}$, by performing a shift of coordinates each time. We furthermore require that $\boldsymbol{0}$ is always an element of $\mathcal{Z}$, i.e. $\mathcal{Z}=\{\boldsymbol{0},\boldsymbol{z}_1,\cdots,\boldsymbol{z}_{q-1}\}$.
	
	Introducing the interpolation set Vandermonde matrix, ${\boldsymbol{\mathrm{Z}}}$, matrix $\boldsymbol{\Pi}^\top$ can be written as
	\begin{equation}
	{\boldsymbol{\Pi}}^\top = \begin{bmatrix}
	\boldsymbol{0} & \boldsymbol{z}_1 & \cdots & \boldsymbol{z}_{q-1} \\
	1 & 1 & \cdots & 1
	\end{bmatrix} = \begin{bmatrix}
	\boldsymbol{0} & \boldsymbol{\mathrm{Z}} \\ 1 & \boldsymbol{1}^\top \\
	\end{bmatrix}
	\end{equation}
	
	As mentioned, in \cite{Conn2008}, the conditions \cref{eq:fullylinear:a} were related to the interpolation set geometry and more specifically to the matrix $\boldsymbol{\mathrm{Z}}$. 
	Here, a generalised lemma is presented establishing similar Taylor-like conditions for \textit{better-than-linear} models, i.e. when $q\geq n+2$. A derivation is included in \cref{appendix:btl}.
	\begin{lemma}
		\label{lemma:bthanlinear}
		Suppose that $m$ is an RBF model with linear tail and RBF contribution $\boldsymbol{\alpha}^\top\boldsymbol{\phi}$ that satisfies the interpolation condition, $m(\boldsymbol{z}_i)=f(\boldsymbol{z}_i)$, for the interpolation set $\mathcal{Z}:=\{\boldsymbol{0},\boldsymbol{z}_1,\dots,\boldsymbol{z}_{q-1}\}\subset\mathcal{B}(\boldsymbol{0},\Delta)$ where $q\geq n+2$ and that satisfies $\|{\boldsymbol{\mathrm{Z}}}\|^{\text{-}1}\Delta\leq\Lambda_\mathrm{Z}$. Furthermore suppose that $f$ and $m$ are continuously differentiably in the hypersphere $\mathcal{B}(\boldsymbol{0},\Delta)$ and that $\nabla f$ and $\nabla\boldsymbol{\phi}^\top \boldsymbol{\alpha}$ are Lipshitz continuous in $\mathcal{B}(\boldsymbol{0},\Delta)$ with Lipshitz constants $\gamma_f$ and $\gamma_\phi$ respectively, then the following two inequalities are satisfied for any $\boldsymbol{x}\in\mathcal{B}(\boldsymbol{0},\Delta)$:
		\begin{subequations}
			\begin{eqnarray}
			|m(\boldsymbol{x})-f(\boldsymbol{x})| &\leq & \left(\tfrac{5}{2}\Lambda_\mathrm{Z}\sqrt{q-1}+\tfrac{1}{2}\right)(\gamma_f+\gamma_\phi)\Delta^2\\
			\|\nabla m(\boldsymbol{x})-\nabla f(\boldsymbol{x})\| &\leq &\tfrac{5}{2}\Lambda_\mathrm{Z}\sqrt{q-1}(\gamma_f+\gamma_\phi)\Delta
			\end{eqnarray}
		\end{subequations}
	\end{lemma}
	
	These conditions are similar to \cref{eq:fullylinear:a}. The Lemma provides therefore an admittedly conservative estimate of the constants $\kappa_f$ and $\kappa_g$, but foremost it reveals the crucial role of the interpolation set geometry, which is reflected through the value of $\|\boldsymbol{\mathrm{Z}}\|$. 
	
	It follows that our interpolation set management must ensure boundedness of $\Delta_k^{\text{-}1}\|\boldsymbol{\mathrm{Z}}_k\|$ from below by the a priori fixed value $\Lambda_\mathrm{Z}^{\text{-}1}$, at each iteration.

	\subsection{Optimal critical set geometry}
		\begin{figure}[t!]
			\centering
			\subfloat[$q=n+1$]{\label{fig:setgeometries:a}
				\includegraphics[trim=1.25cm .9cm 1.5cm 0cm,clip=true,width=.15\textwidth]{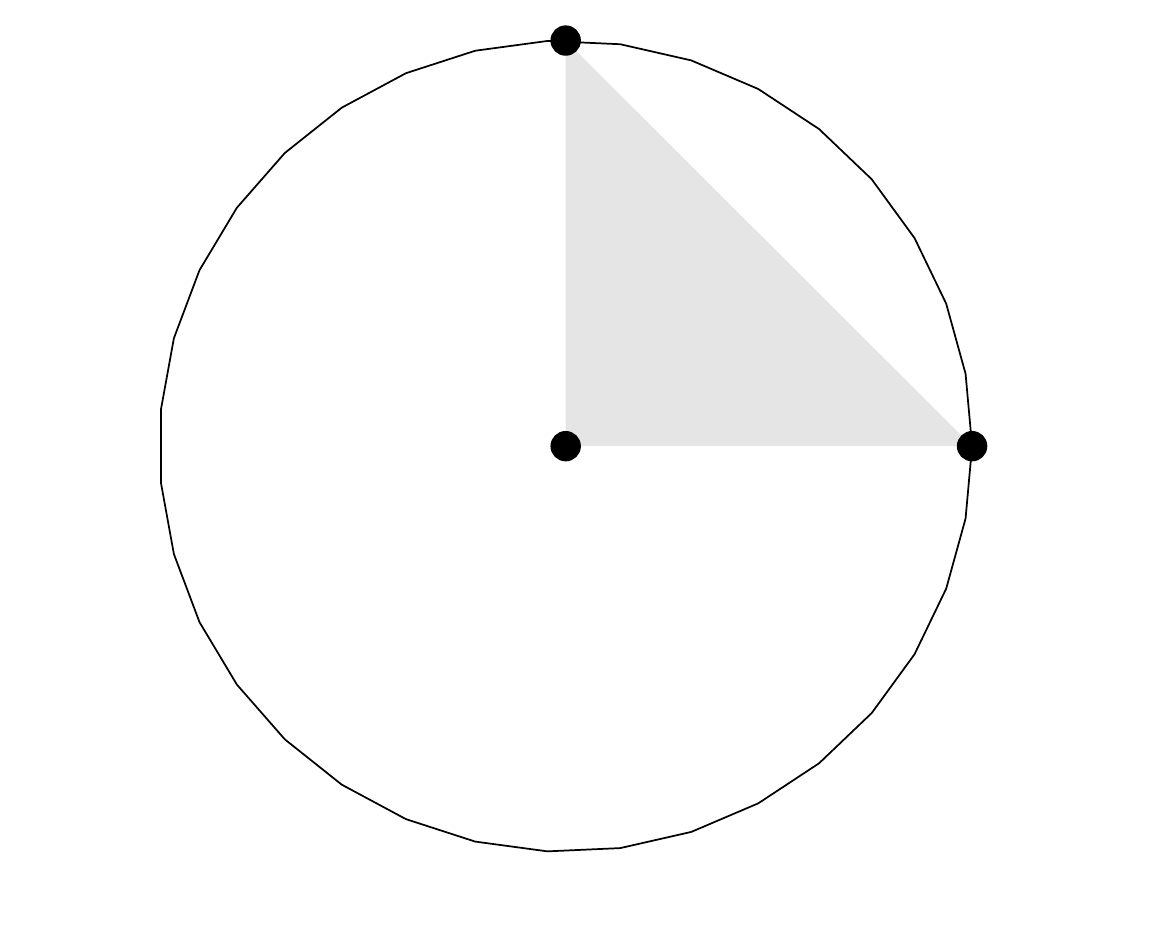}
				\includegraphics[trim=1.25cm .9cm 1.5cm 0cm,clip=true,width=.15\textwidth]{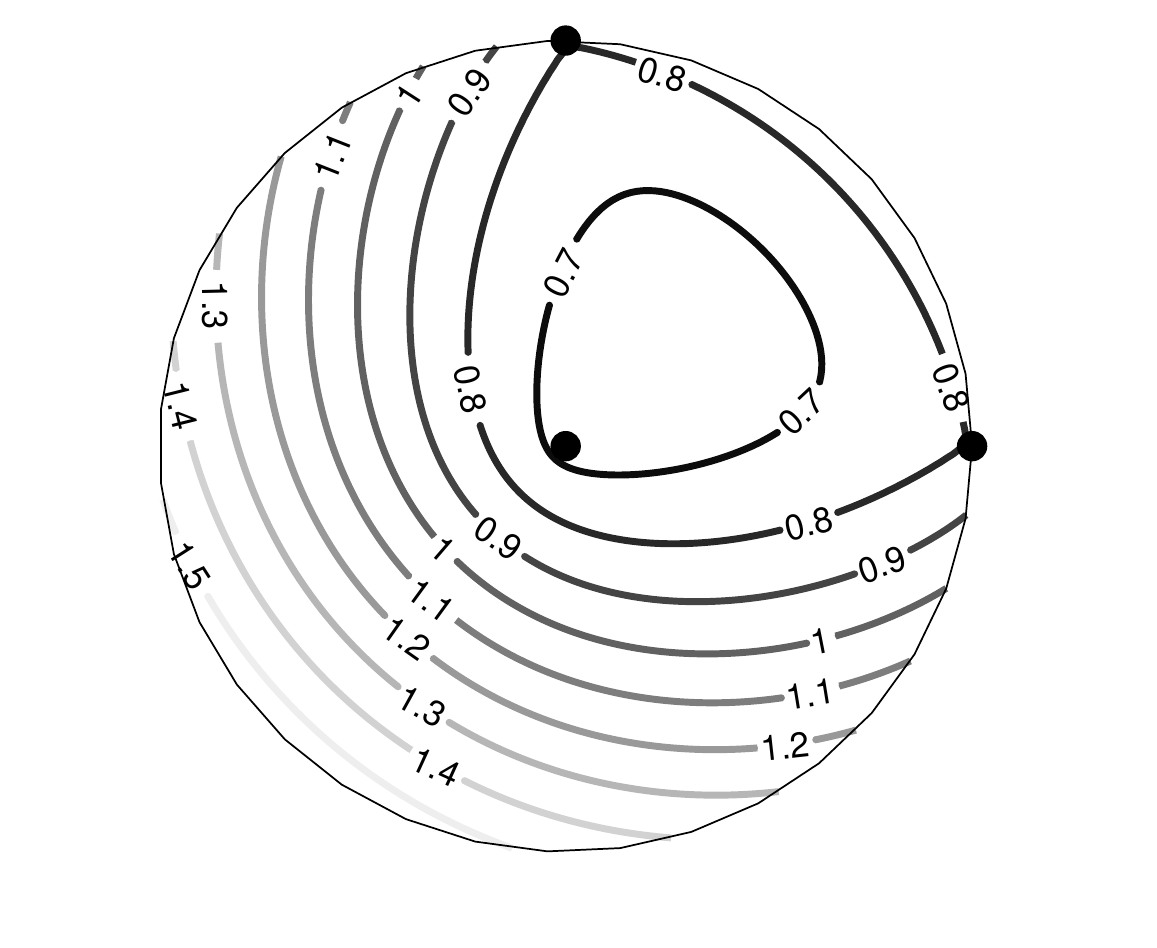}
			}
			\subfloat[$q=n+2$]{\label{fig:setgeometries:b}
				\includegraphics[trim=1.25cm .9cm 1.5cm 0cm,clip=true,width=.15\textwidth]{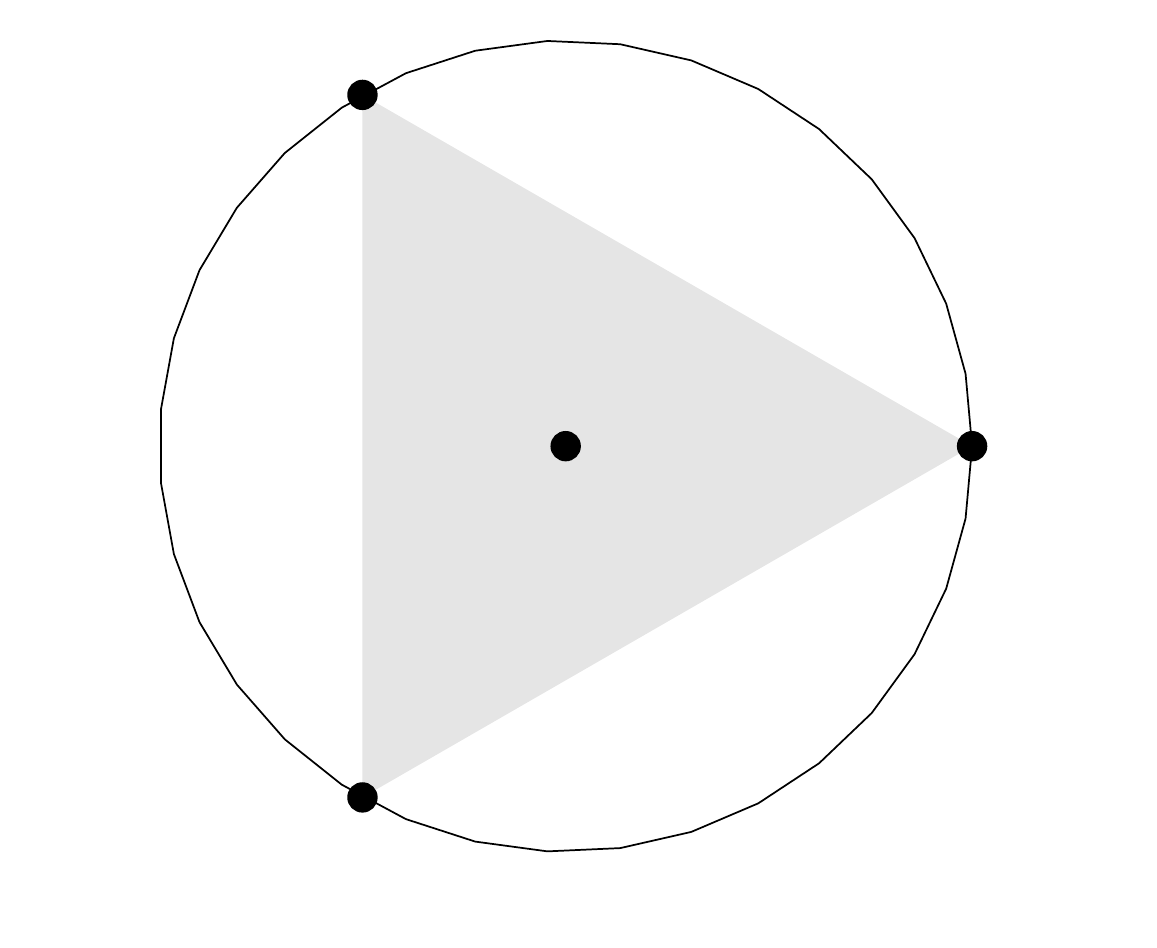}
				\includegraphics[trim=1.25cm .9cm 1.5cm 0cm,clip=true,width=.15\textwidth]{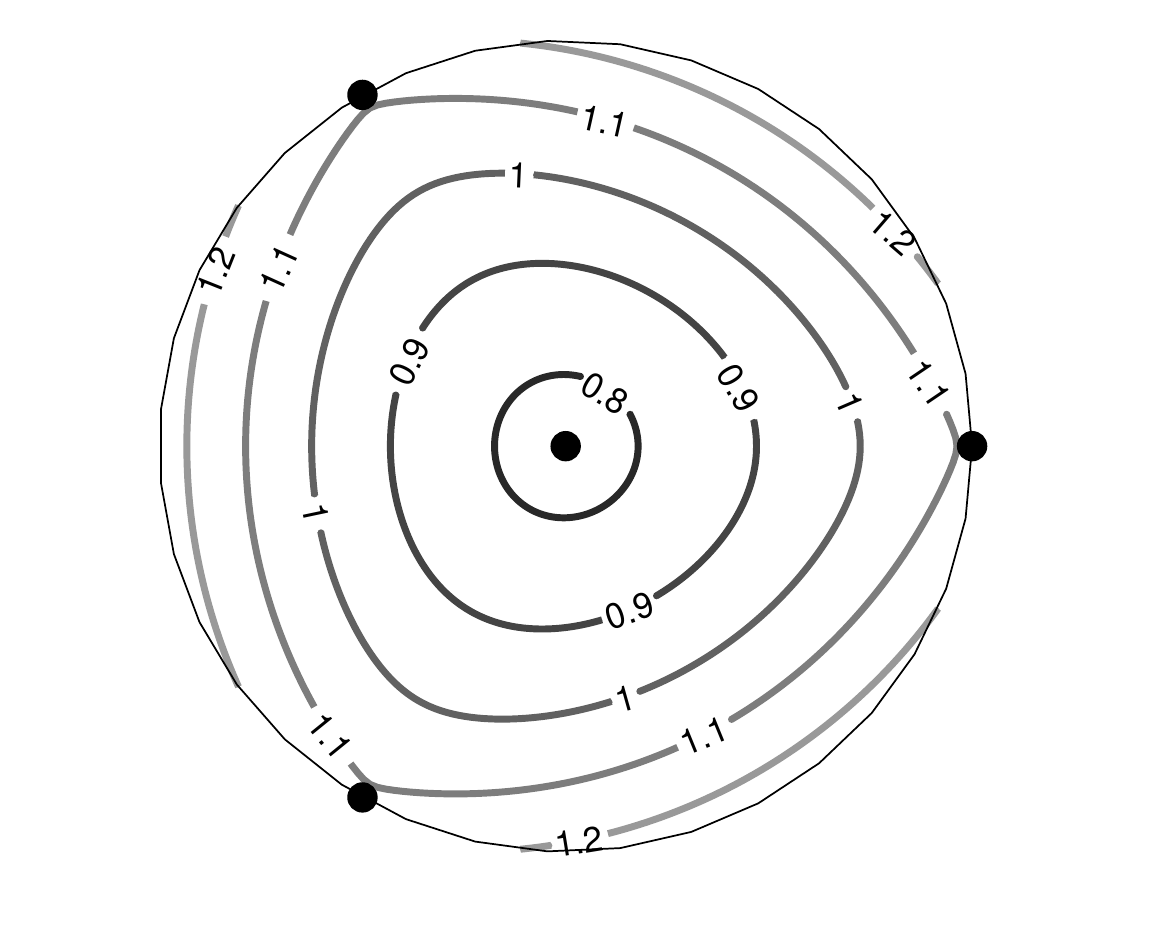}
			}
			\subfloat[$q=2n+1$]{\label{fig:setgeometries:c}
				\includegraphics[trim=1.25cm .9cm 1.5cm 0cm,clip=true,width=.15\textwidth]{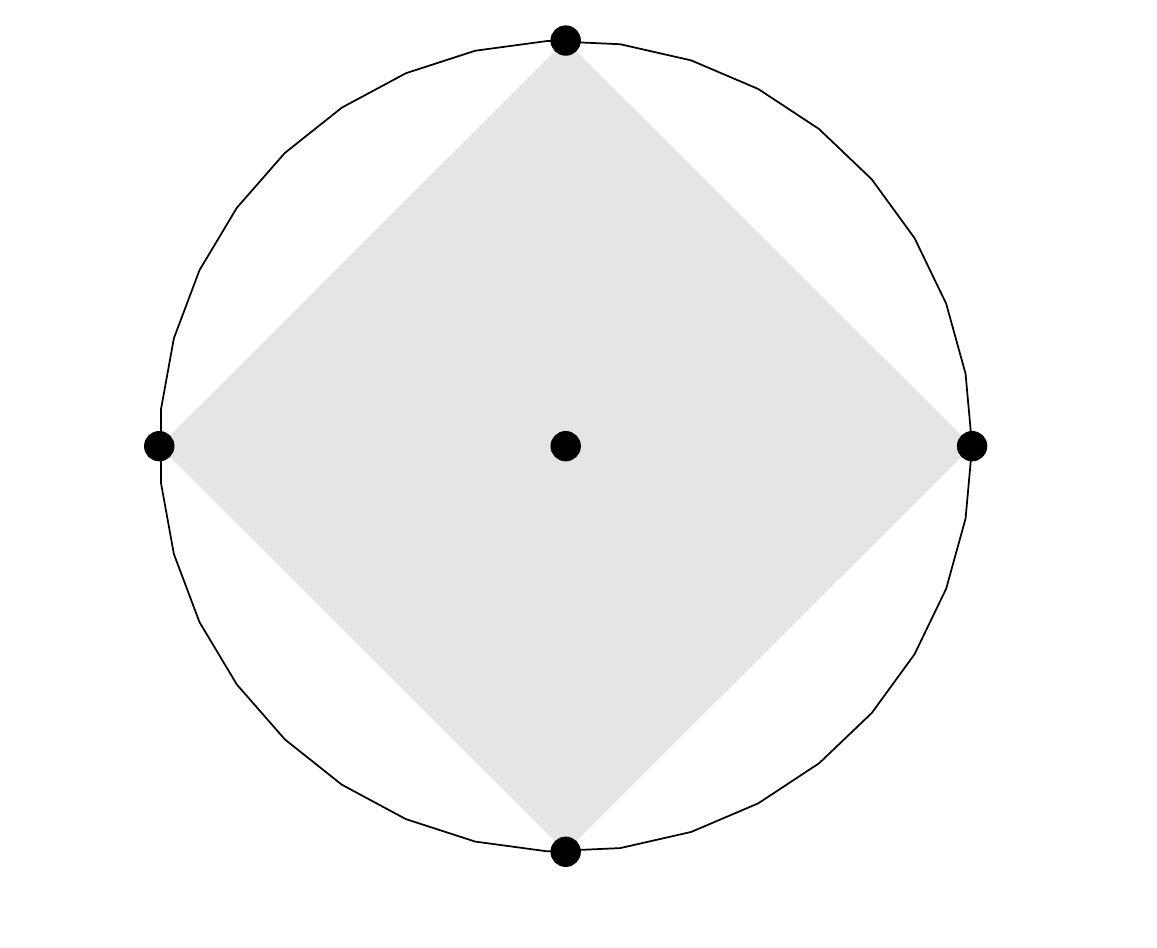}
				\includegraphics[trim=1.25cm .9cm 1.5cm 0cm,clip=true,width=.15\textwidth]{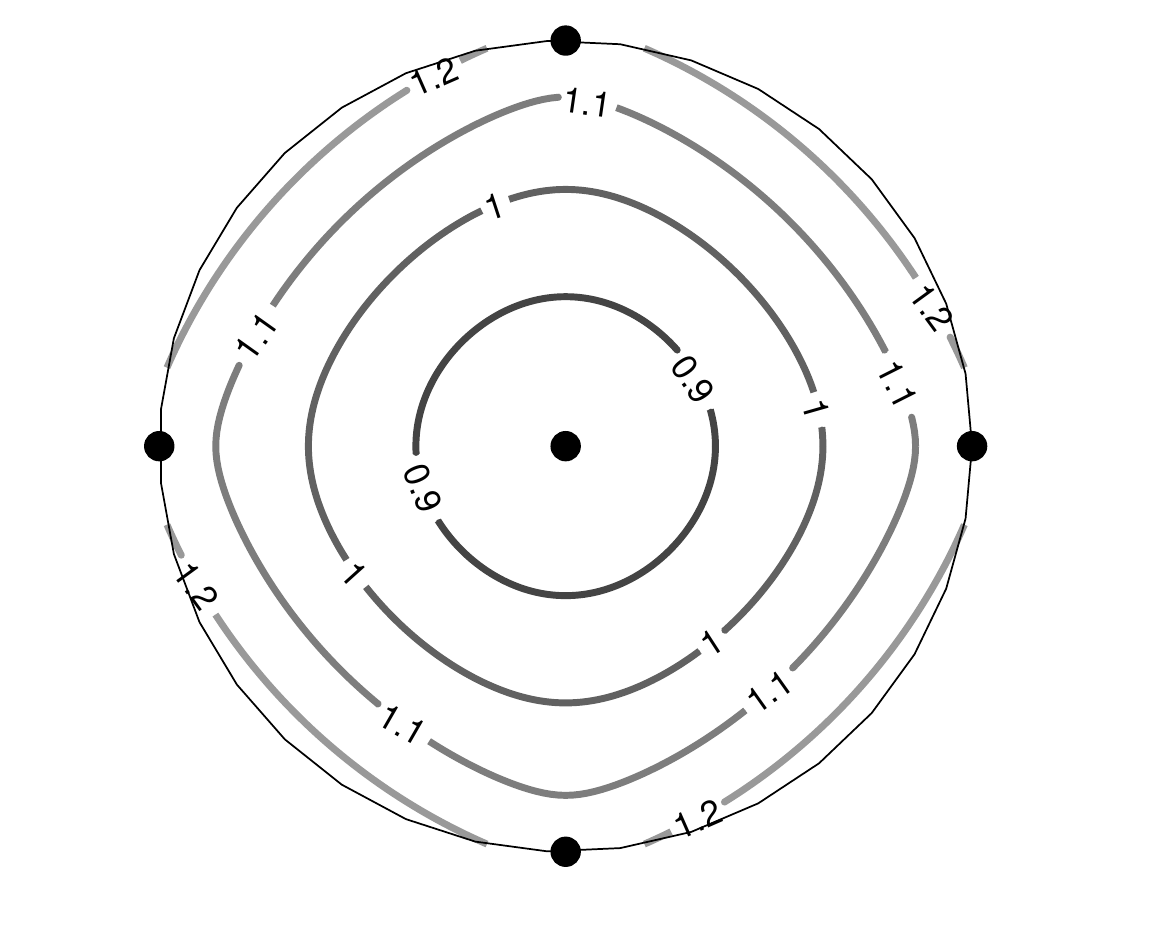}
			}
			\caption{Optimally poised interpolation sets for varying $q$ and $n=2$ next to a contour of $\vartheta(\boldsymbol{x};\mathcal{Z})$. Note that the centre point, $\boldsymbol{0}$, is each time incorporated in the interpolation sets. 
				}
			\label{fig:setgeometries}
			\vspace*{-.5cm}
		\end{figure}
		
	\label{sec:sub:optimal-geometry}
	\label{sec:optimalsetgeom}
	The value of $\Delta^{\text{-}1}\|\boldsymbol{\mathrm{Z}}\|$ is understood as a measure for the algebraic well-poisedness of a given set $\mathcal{Z}\subset\mathcal{B}(\boldsymbol{0},\Delta)$ within $\mathcal{B}(\boldsymbol{0},\Delta)$ \cite{Conn2008}. For any given ammount of interpolation points $q$ we can define the \textit{optimal critical set geometry} that corresponds with the largest value of $\|\boldsymbol{\mathrm{Z}}\|$ corresponding the tightest bounds on $\kappa_f$ and $\kappa_g$. Assuring set well-poisedness hence reduces to generating an affine subset that resembles the optimal set geometry as closely as possible whilst recycling as many points as possible. This strategy balances the requirement to minimize the total amount of function evaluations without corrupting the model quality.
	
	\cref{fig:setgeometries} depicts optimal geometries for $n=2$ dimensions, maximizing the measure $\|\boldsymbol{\mathrm{Z}}\|$ for affine subsets of size $q = n+1$, $q=n+2$ and $q=2n+1$, respectively. The optimal set geometry for $q=n+1$ complies with the normalized set $\Delta^{\text{-}1}\mathcal{Z}\setminus\{\boldsymbol{0}\}$ being orthonormal with respect to the origin, or equivalently when $\mathcal{Z}$ coincides with a standard simplex. This set geometry is pursued by QR pivoting strategies as elaborated in \cite{wild2008orbit}, and Lagrange or Newton polynomials based strategies as elaborated in \cite{oeuvray2009boosters,Conn2008}. The optimal set geometry for $q = n+2$ corresponds with the regular simplex geometry (plus the origin) \cite{nelder1965simplex,gerber1975} and that for $q = 2n+1$ with that of an orthoplex (plus the origin) \cite{crombecq2011novel}. The orthoplex geometry has been recommended by Powell in \cite{powell2006newuoa}. This geometry also determines the fixed search directions in many pattern search optimization methods. 
	
	In this context it is relevant to inspect the spatial distribution of the $2$D geometries in \cref{fig:setgeometries}. The optimal set geometries for $q=n+2$ and $q=2n+1$ are unoriented: the centre of mass of the corresponding sets coincide with the origin. Such does not apply for $q=n+1$ where the set is biased. Note that the above is valid for arbitrary $n$. To quantify this more rigorously, we introduce the metric, $\vartheta^*(\mathcal{Z})$, defined as
	
		\begin{equation}
		\vspace*{-10pt}
		\label{eq:setmetric_v1}
		\vartheta^*\left(\mathcal{Z}\right) = \max_{\boldsymbol{x}\in\mathcal{B}(\boldsymbol{0},\Delta)} \vartheta(\boldsymbol{x};\mathcal{Z}) = \frac{1}{q} \sum_{\boldsymbol{z}_i\in\mathcal{Z}}\|\boldsymbol{x}-\boldsymbol{z}_i\|
		\end{equation}
	
	Its value is equal to the maximum mean distance that any test point in the trust-region can be removed from the modelling set $\mathcal{Z}$ (including $\boldsymbol{0}$). The value gives an indication for the worst-case contribution that the radial basis term (that is directly related to the distance between the test point and the interpolation points) has compared to the linear tail. In \cref{tab:setmetric1}, we compare $\vartheta^*(\mathcal{Z})$ values for different optimal set geometries. It can be observed that $\vartheta^*(\mathcal{Z})$ values can be improved by only adding a single point, referring to the shift from the standard to the regular simplex geometry. Contours of $\vartheta(\boldsymbol{x};\mathcal{Z})$ are presented in \cref{fig:setgeometries} that allow to visualize which regions are not represented well by the optimal set geometries. For $n\rightarrow\infty$, $\vartheta^*(\mathcal{Z})$ tends to $\sqrt{2}$ for all three given configurations, implying that all are as good, or as worse, for infinite dimensional problems. 
	
	The discussion based on $\vartheta^*(\mathcal{Z})$ provides an additional argument to choose the regular simplex over the configuration in \cref{fig:setgeometries:a}, certainly for low dimensional problems. 
	
	\begin{table}[t!]
		\caption{$\vartheta(\mathcal{Z})$ for different optimal set geometries contained within $\mathcal{B}(\boldsymbol{0},1)$.}
		\label{tab:setmetric1}
		\centering
		\begin{tabular}{|c c c c|} \hline
			$q$ & $n+1$ & $n+2$ & $2n+2$ \\ \hline\hline
			$n\geq2$ & $\frac{1+\sqrt{2}\sqrt{n^2+n\sqrt{n}}}{n+1}$ & $\frac{3+\sqrt{2}\sqrt{n^2-n}}{n+2}$ & $\frac{1+\sqrt{2}\sqrt{n^2-n\sqrt{n}}+\sqrt{2}\sqrt{n^2+n\sqrt{n}}}{2n+1}$ \\ \hline\hline
			$2$ & $1.5652$ & $1.25$ & $1.0453$ \\
			$3$ & $1.5821$ & $1.2928$ & $1.1552$  \\
			$4$ & $1.5856$ & $1.3165$ & $1.2142$  \\
			$5$ & $1.5844$ & $1.3321$ & $1.2513$  \\
			$\infty$ & $\sqrt{2}$ & $\sqrt{2}$ & $\sqrt{2}$ \\ \hline
		\end{tabular}
	\end{table}
	
	
	\subsection{Regular simplex affine subsets}
	\label{sec:sub:regular-simplex}
	According to the previous subsection the optimal affine subset (without the origin) for $q=n+2$ is a regular simplex. $\boldsymbol{0}$ remains an element of the affine subset which affects the definition of a regular simplex affine subset based on the $n+1$ remaining interpolation points.
	
	Formally, the set of all regular simplices lying on the surface $\partial\mathcal{B}(\boldsymbol{0},1)$ is given by \cref{eq:regularsets}. $l_n(\cdot)$ (see \cref{appendix:simplex}) is a scaling factor depending on the radius of the $n$-sphere (here $1$).
	\begin{equation}
	\label{eq:regularsets}
	\mathcal{R}^n(\boldsymbol{0},1) = 
	\left\{
	\{\boldsymbol{s}_0,\dots,\boldsymbol{s}_n\}
	\left|
	~\boldsymbol{s}_i\in\partial\mathcal{B}(\boldsymbol{0},1) \text{ and }
	\|\boldsymbol{s}_i - \boldsymbol{s}_j\| = l_n(1), i\neq j
	\right.
	\right\}
	\end{equation}
	
	The geometry of the regular simplex can be related to the matrix norm, $\|\boldsymbol{\mathrm{Z}}\|$, by examining the hypervolume spanned by the convex hull of the affine subset, $\mathcal{Z}\setminus\{\boldsymbol{0}\}$. Consider therefore the hypervolume, $V(\mathcal{Z})$, of an arbitrary set, $\mathcal{Z}=\{\boldsymbol{z}_1,\cdots,\boldsymbol{z}_{n+1}\}$, defined in \cref{eq:hypervolume} with spanning vectors, $\{\boldsymbol{z}_{2}-\boldsymbol{z}_1 ,\cdots,\boldsymbol{z}_{n+1}-\boldsymbol{z}_1\}$. As outlined in \cite{gerber1975,tanner1974some}, the regular simplex has volumetric extremity properties and maximizes the set well-poisedness metric, $\|\boldsymbol{\mathrm{Z}}\|$.
	\begin{equation}
	\label{eq:hypervolume}
	V(\mathcal{Z}) = \frac{\left|\det\left(\begin{bmatrix}
		\boldsymbol{z}_2-\boldsymbol{z}_1 & \cdots & \boldsymbol{z}_{n+1}-\boldsymbol{z}_1
		\end{bmatrix}\right)\right|}{n!}
	\end{equation}
	
	The algebraic problem of maximizing $\|\boldsymbol{\mathrm{Z}}\|$ hence coincides with the geometric problem of maximizing the volume of the convex hull of $\mathcal{Z}$. Approaching the problem from this point of view helps to devise an intuitive interpolation set management resulting in a critical subset that is (close to) an element of the set $\mathcal{R}^n(\boldsymbol{x}_k,\Delta_k)$; whilst recycling interpolation points evaluated in previous iterations.	
	
	Based on the the hypervolume $V(\mathcal{Z})$ of the subset $\mathcal{Z}$ a metric can be construed whether $\mathcal{Z}$ matches a regular simplex sufficiently. The principle idea is that $V(\mathcal{Z})$  should be close to the theoretical maximal hypervolume of a regular simplex, $V^*$, and this within a predefined relative error factor $0\leq \mu_1\leq 1$. If $V(\mathcal{Z})\geq \mu_1 V^*$ is not satisfied, new interpolation points need to be generated. This procedure guarantees $\Delta^{\text{-}1}\|\mathbf{Z}\|\geq\Lambda_\mathrm{Z}^{\text{-}1}$.
	
	
	Before we elaborate the particularities of our set management, we mention that each set of points $\mathcal{R}\in\mathcal{R}^n_m(\boldsymbol{0},1)$, is a lower order regular simplex itself where $m$ is the dimension of the regular subsimplex considered and therefore $m\leq n$. $\mathcal{R}^n_m(\boldsymbol{0},1)$ is defined \cref{eq:subregularsimset}. This statement holds since any subset of a set $\mathcal{R}\in\mathcal{R}^n(\boldsymbol{0},1)$ simply inherits the property that the distance between every two points is the same. 
	\begin{equation}
	\label{eq:subregularsimset}
	\mathcal{R}^n_m(\boldsymbol{0},1) = \left\{\{\boldsymbol{s}_0,\dots,\boldsymbol{s}_m\}\left|\exists\{\boldsymbol{s}_{m+1},\dots,\boldsymbol{s}_n\}:\{\boldsymbol{s}_0,\dots,\boldsymbol{s}_n\}\in\mathcal{R}^n(\boldsymbol{0},1)\right.\right\}
	\end{equation}

	\subsubsection{$\mathbf{Expand2volume}$}
	Based on the previous idea, \cref{algo:e2v}: $\mathrm{Expand2volume}$, constructs an affine subset in a systematic manner so that $m_k$ classifies as a fully linear model for given $\Lambda_Z$. The goal is to attain an affine subset with hypervolume $V(\mathcal{Z})\geq\mu_1 V^*$ whilst adding as few points as possible to the available set of interpolation points. The input is the set of all interpolations points available, scaled with respect to the  trust radius, $\Delta_k$, and unbiased with respect to the center, $\boldsymbol{x}_k$, i.e. $\mathcal{Z}_k = \{ \boldsymbol{z}_i = \Delta_k^{\text{-}1} (\boldsymbol{s}_i - \boldsymbol{x}_k) | \boldsymbol{s}_i \in\mathcal{S}_k\setminus\{\boldsymbol{x}_k\}\}$.  $\mathcal{S}_k$ is the set incorporating every point evaluated so far. The current origin is excluded such that we remain with only the candidate interpolation points to form the affine subset. 
	
	The algorithm embarks on removing all points that are outside a certain periphery, parameterized by the variable $\theta_1\geq1$. Note that we artificially enlarge the trust-region since we allow parameter $\theta_1$ to be greater than 1. (We can always assure that the interpolation model is fully linear by redefining $\kappa_f$ and $\kappa_g$). This ensures that in practice, the number of interpolation points actually accepted to be part of the affine subset is slightly increased, improving the overall convergence. Otherwise, the algorithm might end up expanding an empty set too regularly (certainly for larger $n$). If this operation yields a subset, $\mathcal{Z}_1$, that is still larger than the required affine subset size of $n+1$ (which in practice will almost never occur), all points that are furthest away from the trust surface are additionally removed. 
	If the size of the resulting subset, $\mathcal{Z}_1$, is smaller than $n+1$, it is expanded by calling $\mathrm{Simexpand}$ (\cref{algo:se}).
	
	The procedure described by $\mathrm{Simexpand}$, exploits the idea that if the set $\mathcal{Z}_1$ were to be an element of $\mathcal{R}^n_{n-m}(\boldsymbol{0},1)$, we can expand it with an element of $\mathcal{R}^n_{m-1}(\boldsymbol{0},1)$ so that their union would then constitute an element from $\mathcal{R}^n_{n}(\boldsymbol{0},1)$. The expansion element can be generated by rotating and translating a properly sized lower order regular simplex.
	
	Our numerical implementation proceeds as follows. A set, $\mathcal{X}$, is determined to be an element of $\mathcal{R}^{m-1}(\boldsymbol{0},1)$ embedded in the $n$-dimensional workspace. The direction $\mathbf{e}$ is determined as the orthogonal projection direction of the origin on the subspace spanned by the vertices of $\mathcal{Z}_1$. Subsequently, a rotation matrix $\mathbf{Q}$ is determined so that the basis spanned by the rotated element $\mathbf{Q}\mathcal{X}$ is orthogonal to the union of the basis spanned by $\mathcal{Z}_1$ and the direction $\mathbf{e}$. In a final step the rotated element $\mathbf{Q}\mathcal{X}$ is scaled according to, $r_m^n(1)$, and translated along the direction $\mathbf{e}$ over a distance $d_m^n(1)$. The factors $r_m^n(\cdot)$ and $d_m^n(\cdot)$ are explained in \cref{appendix:simplex} and will only maximize the volume in the occasional coincidence where $\mathcal{Z}_1\in\mathcal{R}^n_{n-m}(\boldsymbol{0},1)$. Nonetheless, this same reasoning can be applied when $\mathcal{Z}_1\notin\mathcal{R}_{n-m}^n(\boldsymbol{0,1})$. In that case an additional subproblem must be solved to identify the optimal value of $\sigma$ so to maximize the volume of their union. Here $\sigma$ is a scaling factor related to $r_m^n(\cdot)$ and $d_m^n(\cdot)$.
	In the description of the algorithm we make no distinction between a set $\mathcal{X}$ and its matrix column vector representation $\mathbf{X}$. 
	
\begin{figure}[t!]
			\centering
			\subfloat[$m=1$]{\label{fig:expand:a}\includegraphics[trim=4cm 8cm 5cm 8cm,clip=true,width=.33\textwidth]{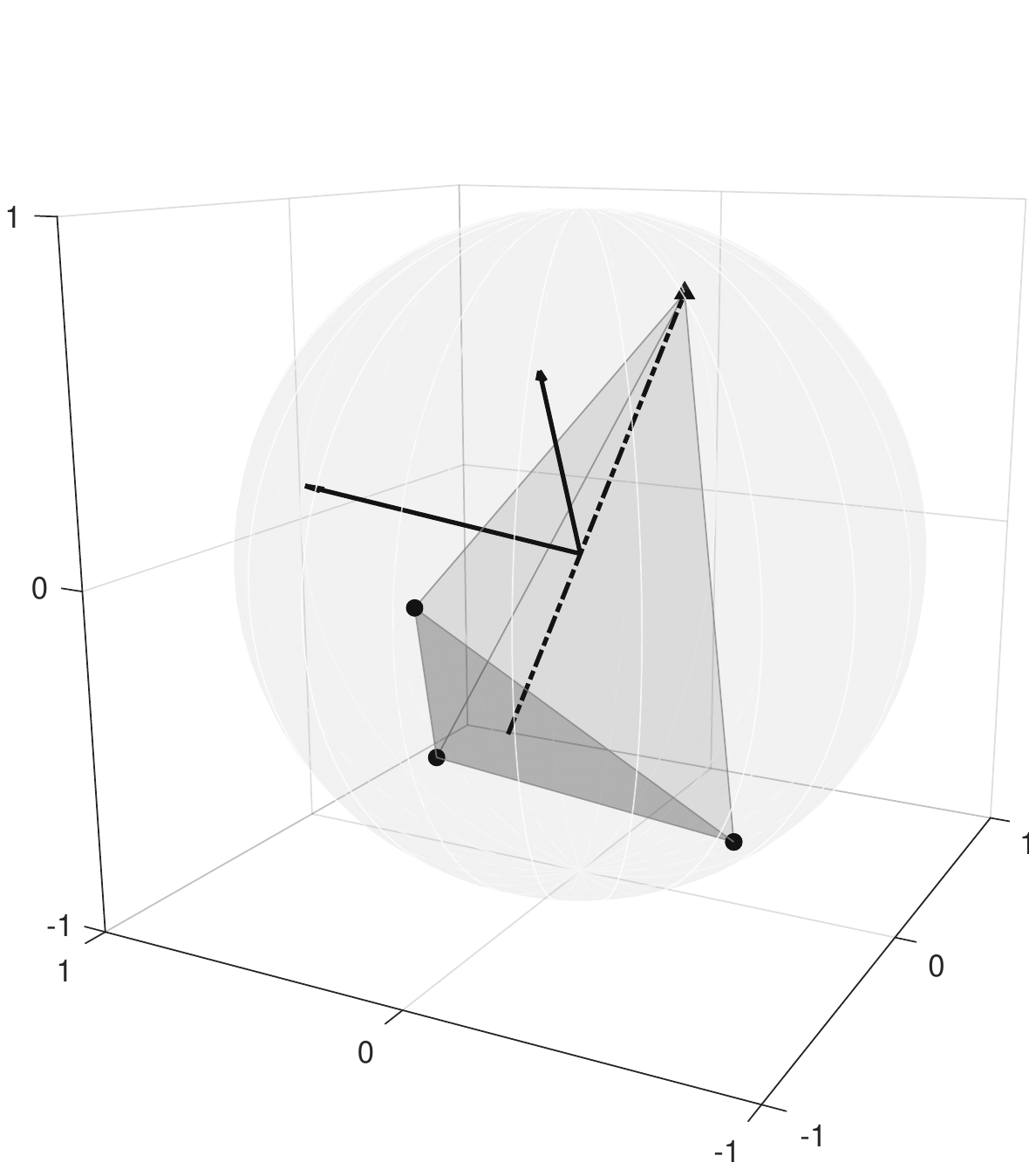}}
			\subfloat[$m=2$]{\label{fig:expand:b}\includegraphics[trim=4cm 8cm 5cm 8cm,clip=true,width=.33\textwidth]{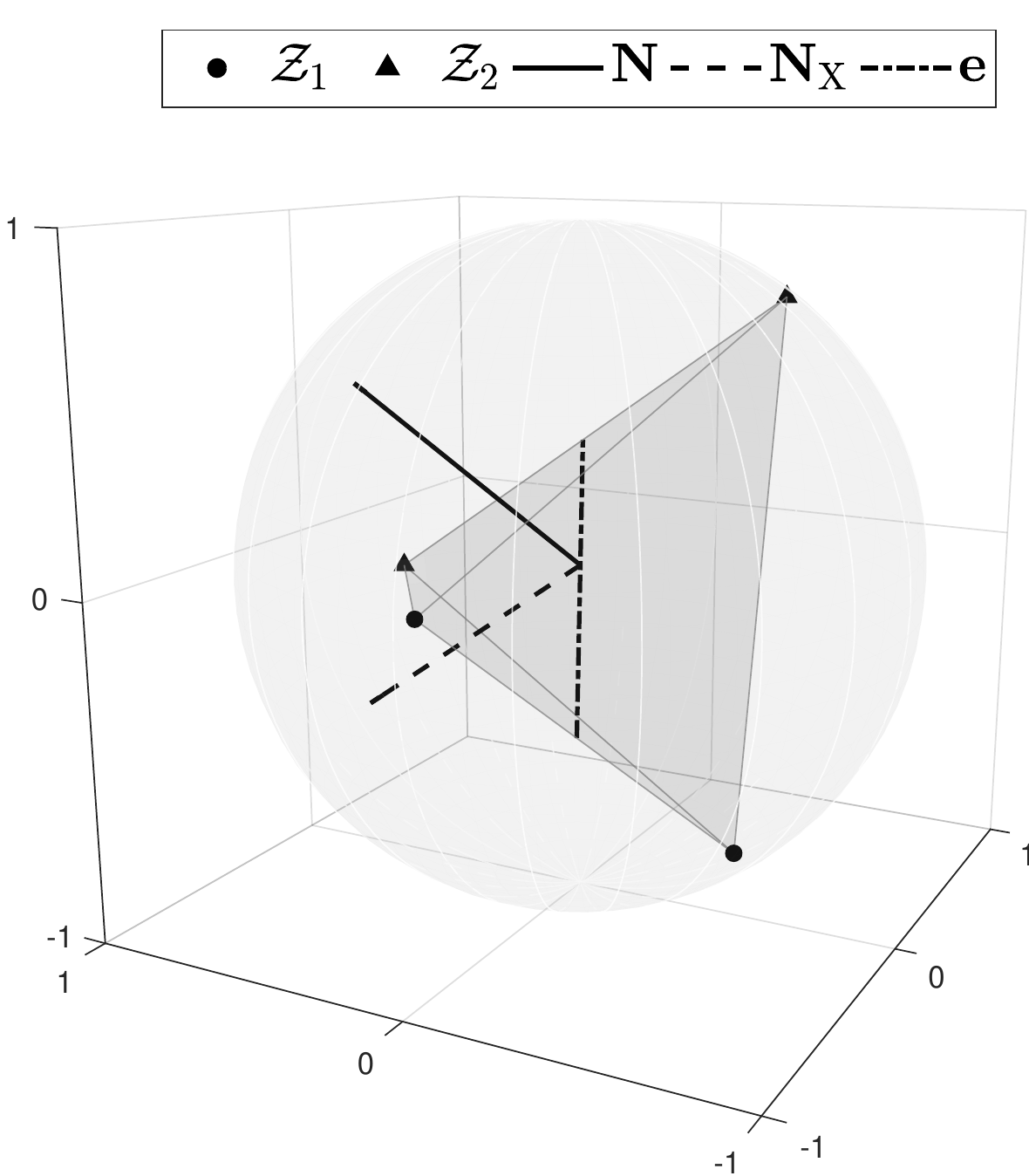}}
			\subfloat[$m=3$]{\label{fig:expand:c}\includegraphics[trim=4cm 8cm 5cm 8cm,clip=true,clip=true,width=.33\textwidth]{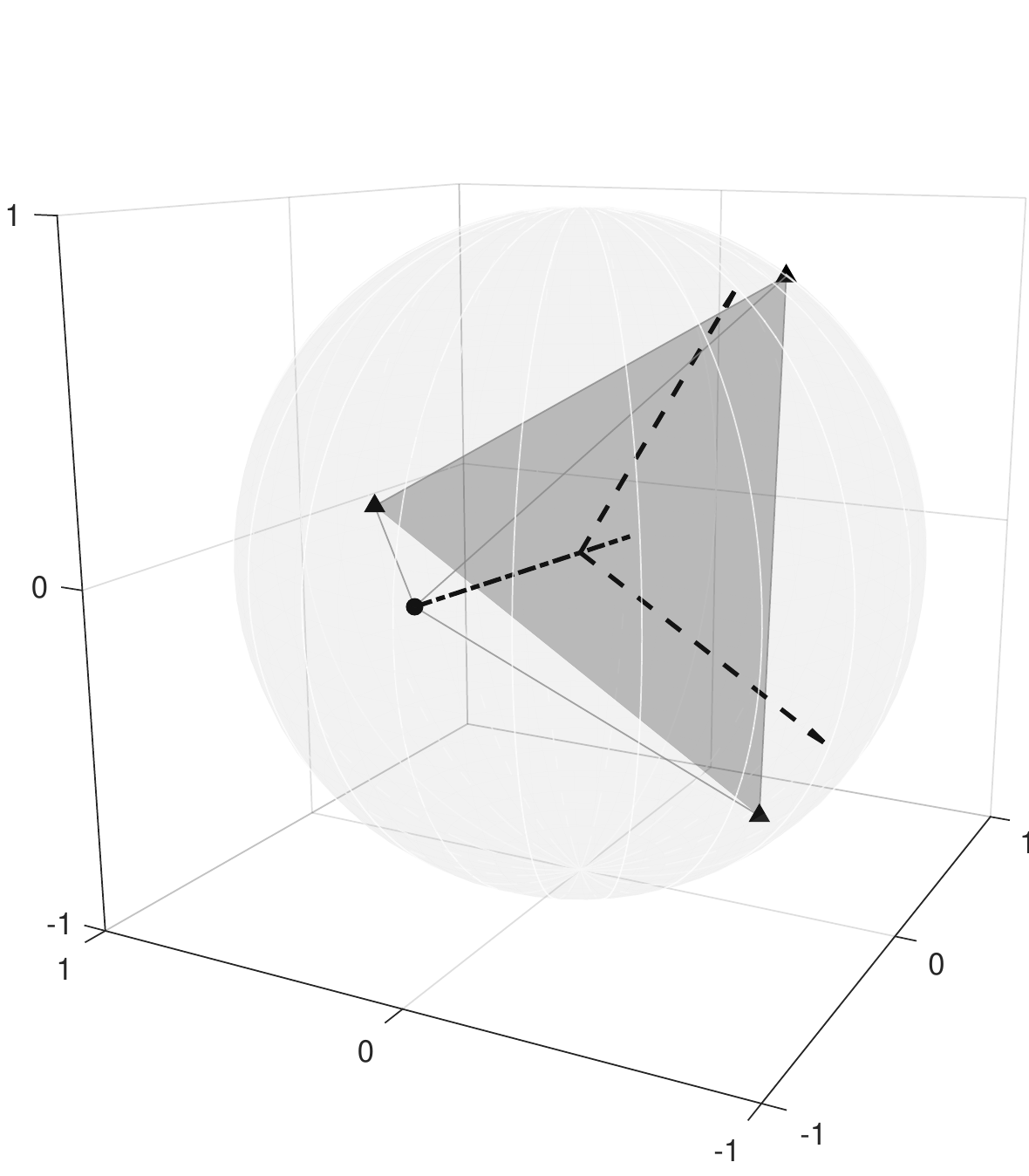}}
			\caption{Illustration of the expansion mechanism of $\mathrm{Simexpand}$ for varying $m$ and $n=3$.}
			\label{fig:expand}
		\end{figure}

	We also use the following operators. Operator, $[\cdot]$, constructs a matrix containing a set of column vectors that span the space determined by the vertices in $\mathcal{Z} = \{\boldsymbol{z}_0,\dots,\boldsymbol{z}_q\}$, i.e. $[\mathcal{Z}] = \begin{bmatrix}
	\boldsymbol{z}_1-\boldsymbol{z}_0& \dots &\boldsymbol{z}_q-\boldsymbol{z}_0
	\end{bmatrix}$.	The function, $\mathrm{regsim}(n,m)$, outputs a set of vertices, $\mathcal{X}  =\{\boldsymbol{x}_0,\dots,\boldsymbol{x}_m\}$, that constitute a regular $m$-simplex in an $n$-dimensional coordinate space. The matrix function, $\mathrm{QR}(\cdot)$, constructs orthonormal bases for the span, ${\boldsymbol{\mathrm{N}}}$, and the kernel, ${\boldsymbol{\mathrm{K}}}$, determined by the columns of the input matrix, respectively. The mechanics of $\mathrm{Simexpand}$ are illustrated in \cref{fig:expand}. Each subfigure depicts a set $\mathcal{Z}_1$ that differs in size from the previous, the corresponding expansion set $\mathcal{Z}_2$, bases $\mathbf{N}$ and $\mathbf{N}_\mathrm{X}$ and the vector $\mathbf{e}$.
	
	If the union $\mathcal{Z}_1\cup\mathcal{Z}_2$ does not satisfy the criticality condition, $V(\mathcal{Z}_1\cup\mathcal{Z}_2)\geq\mu_1V^*$, an iterative procedure is initiated where points are removed from the set $\mathcal{Z}_1$ systematically. The contracted set can then be expanded by applying $\mathrm{Simexpand}$ again. In order to avoid that expansion points are generated too close to a point that was removed in an earlier iteration, information from the removed point is transferred to the contracted set $\mathcal{Z}_1$. Information about a point that gets removed is stored by transferring it to a memory set $\mathcal{Z}'$. The mean of this memory set is used to complement the points in $\mathcal{Z}_1$ that remain, before expanding the set with a new regular $m+1$-simplex. In this fashion, the added simplex is generated away from the mass centre of the original set $\mathcal{Z}_1$.
	
	The algorithm ends after at most $n+1$ iterations, that is when the entire original set $\mathcal{Z}_1$ is transferred to the memory set $\mathcal{Z}'$. We remark that in this case the mass center is expanded with an $(n-1)$-regular simplex that is however not necessarily an element of $\mathcal{R}^n_{n-1}(\boldsymbol{0},1)$. As a consequence, the eventual affine subset needs not to be an element of $\mathcal{R}^n(\boldsymbol{0},1)$. As, this does not ensure that the sufficiency condition is satisfied, the algorithm closes with a final verification of the hypervolume and possible expansion of the set $\mathcal{Z}_2$.
	
	
	\begin{algorithm}[t]
		\caption{$\mathrm{Expand2volume}$}
		\label{algo:e2v}
		\begin{algorithmic}[1]
			\renewcommand{\algorithmicrequire}{\textbf{Input:}}
			\renewcommand{\algorithmicensure}{\textbf{Output:}}
			\REQUIRE $\mathcal{Z}_0 = \mathcal{Z}/\boldsymbol{0},\theta_1,\mu_1$
			\ENSURE $\mathcal{Z}_2$
			\STATE{obtain $\mathcal{Z}_1$ with elements $\boldsymbol{z}_i\in\mathcal{Z}_0$ such that $\|\boldsymbol{z}_i\|\leq\theta_1$}
			\STATE{define $\mathcal{Z}_2 = \emptyset$}
			\IF{$|\mathcal{Z}_1| > n+1$}
			\STATE{order $\boldsymbol{z}_i$ such that $\left|1-\|\boldsymbol{z}_1\|\right|\leq\dots\leq\left|1-\|\boldsymbol{z}_q\|\right|$}
			\STATE{obtain $\mathcal{Z}_1 = \mathcal{Z}_1 / \left\{ \boldsymbol{z}_{n+2},\dots,\boldsymbol{z}_{q} \right\}$}
			\ELSIF{$|\mathcal{Z}_1| < n+1$}
			\STATE{get $\mathcal{Z}_2 = \mathrm{Simexpand}(\mathcal{Z}_1)$}
			\ENDIF
			\STATE{find initial memory $\boldsymbol{z}' = \boldsymbol{z}_i\in\mathcal{Z}_1:\|\boldsymbol{z}_i\|\leq\|\boldsymbol{z}_j\|,\forall j$}
			\STATE{define memory set $\mathcal{Z}' = \boldsymbol{z}'$}
			\WHILE{$V(\mathcal{Z}_1\cup\mathcal{Z}_2)<\mu_1 V^*$ or $|\mathcal{Z}_1|>0$}
			\STATE{find $\boldsymbol{z}_i\in\mathcal{Z}_1:\|\boldsymbol{z}_i-\boldsymbol{z}'\|\leq\|\boldsymbol{z}_j-\boldsymbol{z}'\|,\forall j$}
			\STATE{contract $\mathcal{Z}_1 = \mathcal{Z}_1/\boldsymbol{z}_i$}
			\STATE{redefine memory set $\mathcal{Z}' = \mathcal{Z}'\cup\boldsymbol{z}_i$}
			\STATE{expand $\mathcal{Z}_2 = \mathrm{Simexpand}(\mathcal{Z}_1\cup\mathrm{mean}(\mathcal{Z}'))$}
			\ENDWHILE
			\IF{$V(\mathcal{Z}_1\cup\mathcal{Z}_2)<\mu_1 V^*$}
			\STATE{$\mathcal{Z}_2 = \mathcal{Z}_2\cup\mathrm{Simexpand}(\mathcal{Z}_2)$}
			\ENDIF
			\RETURN $\mathcal{Z}_2$ 
		\end{algorithmic} 
	\end{algorithm}
	
	\begin{algorithm}
		\caption{$\mathrm{Simexpand}$}
		\label{algo:se}
		\begin{algorithmic}[1]
			\renewcommand{\algorithmicrequire}{\textbf{Input:}}
			\renewcommand{\algorithmicensure}{\textbf{Output:}}
			\REQUIRE $\mathcal{Z}_1$
			\ENSURE $\mathcal{Z}_2$
			\STATE{$m = n-|\mathcal{Z}_1|$}
			\STATE{$\mathcal{X} = \mathrm{regsim}(n,m-1)$}
			\STATE{$[{\boldsymbol{\mathrm{N}}},\sim] = \mathrm{QR}([\mathcal{Z}_1])$, $[{\boldsymbol{\mathrm{N}}}_0,{\boldsymbol{\mathrm{K}}}_0] = \mathrm{QR}([\mathcal{Z}_1\cup\boldsymbol{0}])$, $[{\boldsymbol{\mathrm{N}}}_\mathrm{X},{\boldsymbol{\mathrm{K}}}_\mathrm{X}]= \mathrm{QR}([\mathcal{X}\cup\boldsymbol{0}])$}
			\STATE{$[\sim,\mathbf{e}] = \mathrm{QR}([{\boldsymbol{\mathrm{N}}} ~ {\boldsymbol{\mathrm{K}}}_0])$}
			\STATE{${\boldsymbol{\mathrm{Q}}}=[{\boldsymbol{\mathrm{K}}}_0 ~ {\boldsymbol{\mathrm{N}}}_0]\cdot[{\boldsymbol{\mathrm{N}}}_\mathrm{X} ~ {\boldsymbol{\mathrm{K}}}_\mathrm{X}]^\top$}
			\STATE{$\sigma^* = \arg\min\limits_\sigma ~V\left(\sigma\cdot\boldsymbol{1}^\top_m\cdot\boldsymbol{e} + \sqrt{1-\sigma^2}\cdot{\boldsymbol{\mathrm{Q}}} \cdot \mathcal{X}\right)$}
			\STATE{$\mathcal{Z}_2 = [\boldsymbol{z}_1 ~ \dots  ~ \boldsymbol{z}_{n-|\mathcal{Z}|}] =\sigma^* \cdot\boldsymbol{1}^\top_m\cdot\boldsymbol{e} +  \sqrt{1-\sigma^{*2}} \cdot {\boldsymbol{\mathrm{Q}}} \cdot \mathcal{X}$}
			\RETURN $\mathcal{Z}_2$ 
		\end{algorithmic} 
	\end{algorithm}

	\subsection{Remarks}
	With respect to \cref{algo:e2v} with \cref{algo:se} described above, we add the following. During the initialization of $\mathrm{Expand2volume}$, all points that are furthest away from the trust-surface are removed. This is a heuristic. Alternatively, one might run through a more complex scheme where every possible $n+1$ size set is verified with respect to their hypervolume and the largest set is selected. These points could also be added to the memory set. However, this is a situation that will almost never occur in practice, partially because the optimization algorithm will guide the iterations to unexplored regions of the parameter space, also because the trust region will contract when the sequence of iterates starts to converge. So by proper choice of $\mu_1$ and $\theta_1$ it can be avoided that the complete interpolation set of the previous iteration is within the current periphery.
	
	Furthermore, we remark that alignment of the kernel associated to the space spanned by the union $\mathcal{Z}_1\cup\boldsymbol{0}$ and the space spanned by the regular simplex $\mathcal{X}$ is underdetermined. This additional degree of freedom remains here unexploited. Alternatively, this degree of freedom could be employed by orienting the added simplex as such that the generated points are located furthest away from the available set of interpolation points or such that the added points correspond to a mean minimal predicted objective function value. 
	

	\newpage
	\section{Ellipsoidal trust-region framework}
	\label{sec:ellips}
	The procedures described in the previous two sections can be employed to verify and construct a model, ${m}_k$, that satisfies \cref{def:fullylinear} for a fully linear model in a given trust-region, $\mathcal{B}(\boldsymbol{x}_k,\Delta_k)$. Hence, these procedures can be used to embed an interpolation framework in the trust-region algorithm with fully linear models described by \cref{alg:TR:fullylinear}, factually realising a IBO algorithm.
	
	As for now the geometry of the interpolation points is considered only with regard to the input space and with respect to the trust-region defined in it. Given a spherical trust-region, our input space oriented sampling strategy will treat every direction equally in function of a well-conditioned system matrix $\mathbf{Z}$. However, from the perspective of optimization, it might be of interest to account for the geometry of the objective function as well. 
	
	The gradient of the objective function will be more sensitive to certain search directions then to others as determined by the local curvature of the objective function. Currently the set management however treats all directions equally (as a result of the spherical shape of the trust-region) and will try to homogenize the distance from the trust-region center and in between interpolation points. Now as a result of the local curvature, the corresponding function values may vary faster either slower than can be expected when solely considering the local objective gradient. Hence the geometry of the objective function with respect to the arbitrary choice of the optimization space coordinates is reflected by its local curvature.
	
	Therefore we propose to execute the set management in a transformed space by introducing an auxiliary optimization parameter, $\boldsymbol{u}$. The goal of this transformation is to obtain a normalized curvature of the objective function in the auxiliary optimization space. We propose to use an adaptive affine mapping from the auxiliary space to the original space, that changes with the trust-region iterations. Notice that the origin in the auxiliary space corresponds with the trust-region centre in the original space. 
	\begin{equation}
	\label{eq:affinetrans}
	\boldsymbol{x}(\boldsymbol{z}) = \mathbf{A}_k \boldsymbol{u} + \boldsymbol{b}_k = {\boldsymbol{\mathrm{A}}}_k\boldsymbol{u} + \boldsymbol{x}_k 
	\vspace*{-5pt}
	\end{equation}
	
	The trust-region will be stretched along those directions where the local gradient is varying slowly and shrinked in those directions where the local gradient is varying rapidly. As a consequence, objective function values evaluated along the trust-region boundary will increase or decrease in equal amounts with respect to the objective function value at the trust-region centre. Such will have a positive effect on the hyperparameter tuning described in \ref{sec:sub:hyper} since the principle directions are now homogenized.
	
\begin{figure}[b!]
		\centering
		\subfloat[spherical region]{\label{fig:trusregions:a}
			\includegraphics[width=0.4\textwidth]{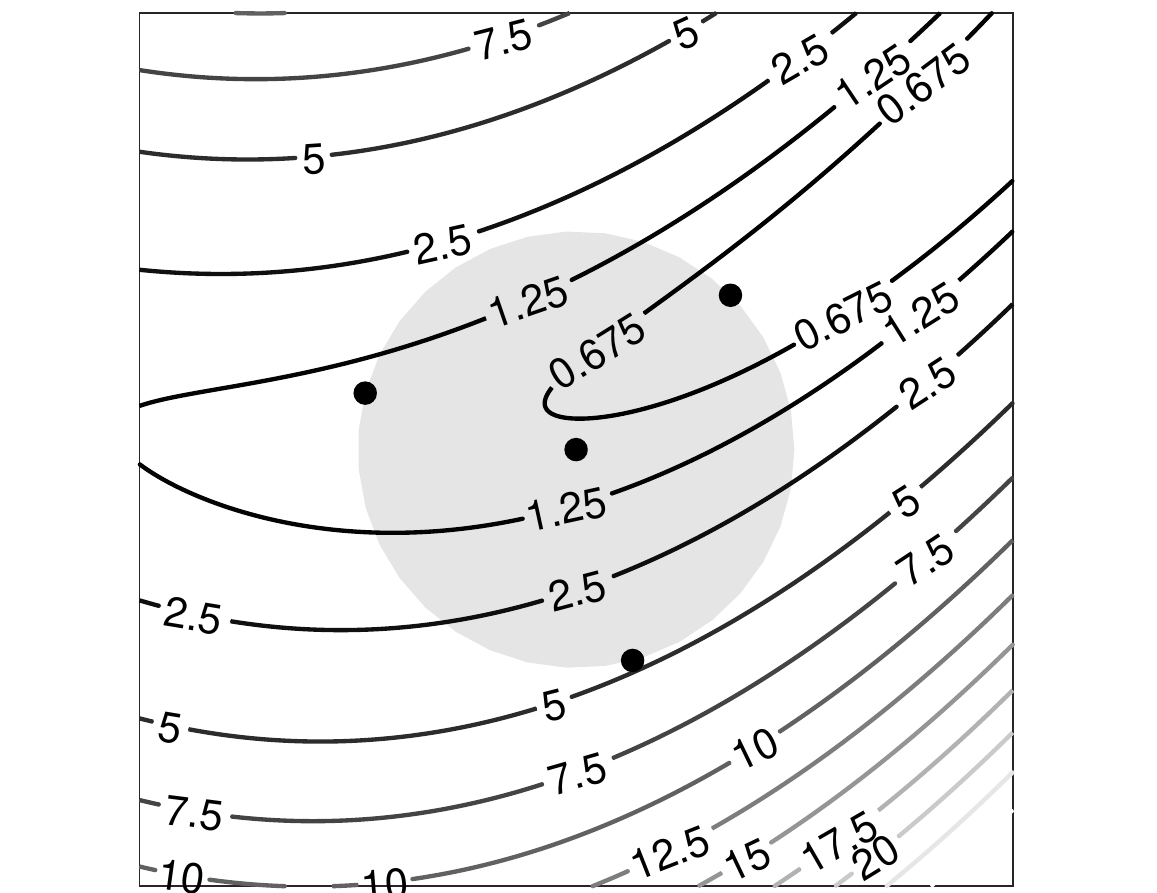}
		}
		\subfloat[normalized ellipsoidal region]{\label{fig:trusregions:b}
			\includegraphics[width=0.4\textwidth]{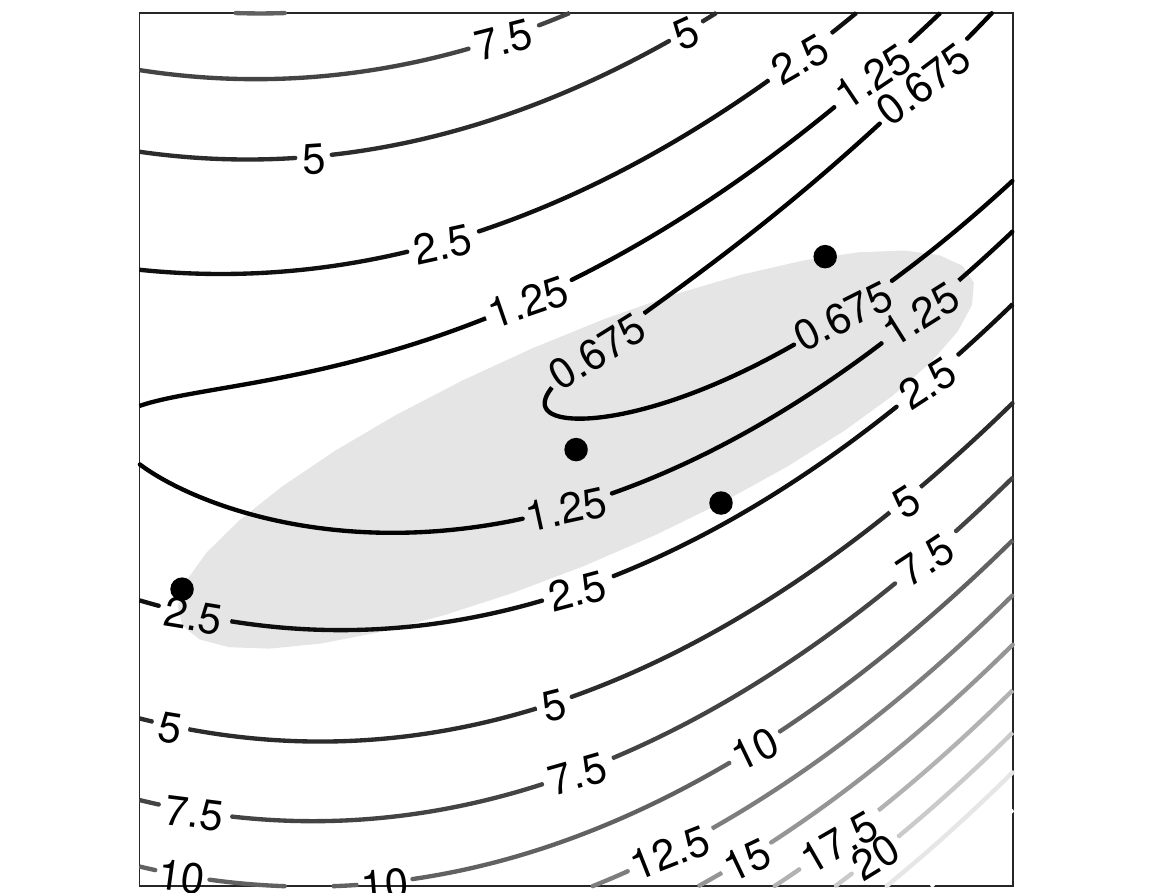}
		}
		\caption{Illustration of a spherical trust-region and its isovolumetric ellipsoidal counterpart for the two dimensional Rosenbrock function. The curvature was normalized using the exact Hessian.}
		\label{fig:trusregions}
		\vspace*{-.5cm}
	\end{figure}

	Another consequence is that a spherical trust-region in the auxiliary space $\boldsymbol{z}$ will now have an ellipsoidal shape in the original optimization space illustrated by \cref{fig:trusregions:a} and \ref{fig:trusregions:b}. As a result, a step of length $1$ taken in the auxiliary space represented in the original space may exceed the original spherical trust-region radius. Two situations can be distinguished when inspecting the directional curvature. When having a large directional curvature on the one hand, the ellipsoidal geometry will have a stabilizing effect considering that the set of admissible steps (i.e. the trust-region) becomes less aggressive in terms of large jumps in the objective function, which will also have an indirect effect on the modelling procedure. On the other hand when having a small directional curvature, the ellipsoidal trust region will allow the optimization procedure to take large steps in the direction of a slowly varying gradient. This way the original optimization space will be explored more rapidly. For instance in \cref{fig:trusregions:b}, the trust-region in the auxiliary space encloses a larger part of the inner contour that contains the actual minimum.
	
	\subsection{Curvature normalization}
	The goal of the coordinate transformation is to normalize the local curvature of the problem in the new coordinate space. This condition can be expressed mathematically as given. For reasons that will become apparent, we introduce an additional scaling factor $\lambda_k^{2}$.
	
	\begin{equation}
	\nabla^2_{\boldsymbol{u}}f(\boldsymbol{x}_k) = {\boldsymbol{\mathrm{A}}}_k^\top \nabla^2_{\boldsymbol{x}}f(\boldsymbol{x}_{k}) {\boldsymbol{\mathrm{A}}}_k = \lambda_k^{2}\mathbf{I}_n
	\end{equation}
	
	Suppose now that we have an estimate ${\boldsymbol{\mathrm{B}}}_{k}$ for the hessian $\nabla^2_{\boldsymbol{x}}f(\boldsymbol{x}_{k})$. We can formulate an expression for ${\boldsymbol{\mathrm{A}}}_k$ by considering the Singular Value Decomposition (SVD) of the matrix ${\boldsymbol{\mathrm{B}}}_{k}$. The SVD decomposition of the positive definite matrix $\mathbf{B}_k$ is given by $\mathbf{U}_k \boldsymbol{\Lambda}_k \mathbf{U}_k^\top$. Here $\mathbf{U}_k$ is an orthonormal matrix and ${\boldsymbol{\Lambda}}_k$ is a diagonal matrix with the singular values $\sigma_1,\cdots,\sigma_n$ on its main diagonal. It is then easily verified that 
	\begin{equation}
	\mathbf{A}_k = \lambda_k \mathbf{U}_k
	\boldsymbol{\Lambda}_k^{\text{-}\frac{1}{2}}
	\end{equation}
	
	Finally we will utilise the scaling factor $\lambda_k$ to keep subproblem $\cref{eq:TR:subproblem}$ as it manifests in the auxiliary coordinate space proportional to its manifestation in the original coordinate space. Mathematically this proportionality can be realised by determining a measure that remains invariant under the transformation. Recall that as a result of the affine coordinate transformation, a spherical trust-region in the auxiliary optimization space will adopt an ellipsoidal geometry in the original optimization space. We propose the following measure of proportionality: the value of $\lambda_k$ is determined such that the spherical and ellipsoidal trust-region representation remain isovolumetric. 
	
	Consider that the volume of an $n$-dimensional hypersphere with radius $\Delta_k$, is equal to $\eta(n)\cdot\Delta_k^n$, and that the volume of an $n$-dimensional ellipsoid with given half principle axis lengths, $a_1,\dots,a_n$, is equal to $\eta(n)\cdot\prod_i a_i$, where $\eta_n$ is a mutual scaling factor. We also have that by construction of $\mathbf{A}_k$, half the principle axis lengths can be written in function of the singular values as $a_i = \lambda_k\sigma_i^{\text{-}1/2}$. An isovolumetric relation between the trust-region representations in both the original as the auxiliary coordinate space can be maintained by proper choice of the scaling factor $\lambda_k$, that is
	\begin{equation}
	\lambda_k = \Delta_k\prod_{i=1}^{n}
	\sigma_i^{\tfrac{1}{2n}} 
	\end{equation}
	
	The exact coordinate transformation of \cref{eq:affinetrans} can then be written explicitly as
	\begin{equation}
	\boldsymbol{x}(\boldsymbol{u}) = \Delta_k\prod_{i=1}^{n}{\sigma_i}^{\frac{1}{2n}}{\boldsymbol{\mathrm{U}}}_k{\boldsymbol{\Lambda}}_k^{\text{-}\frac{1}{2}}\boldsymbol{u} + \boldsymbol{x}_k 
	\end{equation}
	
	\subsection{Filtered Hessian approximation}
	Several methods can be devised to provide an estimate of the Hessian, $\mathbf{B}_k$. Since our sample strategyensures that the local interpolation model is curved, we can calculate an approximation of the Hessian at the next iterate in the auxiliary coordinate space directly from the model, i.e. $\nabla_{\boldsymbol{u}}^2 m_k(\boldsymbol{u}_{k+1})$. 
	
	This estimate can be calculated back to the original coordinate space so to provide an estimate for the actual Hessian, $\hat{\mathbf{B}}_{k+1} = \mathbf{A}_k\nabla_{\boldsymbol{u}}^2 m_k(\boldsymbol{u}_{k+1})\mathbf{A}_k^\top$. In order to get a smooth approximation, a low-pass filter is proposed, parameterised by $\mu_2$. Note that if $\mu_2 = 1$, this framework reduces to the original spherical trust-region framework.
	\begin{equation}
	\label{eq:hessianfilter}
	\mathbf{B}_{k+1} = \mu_2 \mathbf{B}_k + (1-\mu_2)\hat{\mathbf{B}}_{k+1}, ~ \mathbf{B}_0 = \mathbf{I}_n
	\end{equation}
	
	Alternatively, the well-established BFGS-update can be used
	\begin{equation}
	\mathbf{B}_{k+1} = \mathbf{B}_{k} + \frac{\boldsymbol{y}_k \boldsymbol{y}_k^\top}{\boldsymbol{y}_k^\top \boldsymbol{s}_k} - \frac{\mathbf{B}_k\boldsymbol{s}_k\boldsymbol{s}_k^\top\mathbf{B}_k}{\boldsymbol{s}_k^\top \mathbf{B}_k \boldsymbol{s}_k} 
	\end{equation}
	where $\boldsymbol{y}_k = \nabla {m}_k(\boldsymbol{x}_k + \boldsymbol{s}_k) - \nabla {m}_k(\boldsymbol{x}_k)$.
	
	However we argue against the use of this strategy given the current framework. As a result of the method used to construct the surrogate model, large trust-regions can be anticipated which would detoriate the approximation principle of the BFGS update. Moreover, since the gradient calculated to compute $\boldsymbol{y}_k$ are themselves an approximation of the real values, the resulting Hessian approximation becomes unreliable as confirmed by our numerical experiments.

	\section{Numerical results}
	\label{sec:numericalresults}
	In this section, we present numerical results of the presented derivative-free optimization algorithm applied on a set of unconstrained problems from the CUTEst test environment \cite{gould2015cutest}. As indicated by the neighbourhood representation metric, $\vartheta(\mathcal{Z})$, the method should be especially suited for low dimensional problems and was developed from this point of view. We have applied the algorithm on a set of test functions with a maximum problem dimension $n=8$. All tested algorithms were initialized from the traditional starting point \cite{more1981testing}. The following benchmark procedure has been adopted and results are presented in \cref{tab:benchmark}.
	
	The results as obtained with $4$ different versions of the present algorithm have been compared to those obtained with its closest competitors in the context of derivative-free optimization, being ORBIT \cite{wild2008orbit} and NMSMAX. Since ORBIT outperforms the NEWUOA algorithm (at least for small scale problems \cite{wild2008orbit}) and given that no \MATLAB~ implementation is available, we did not compare our algorithm with NEWUOA \cite{powell2006newuoa}. The details of ORBIT are elaborated in the associated paper and are essentially similar to the details of the presented algorithm when the spherical trust-region framework is used. That is apart from the novel regular simplex based sampling strategy. We used the implementation made available at \cite{orbit_mat}. The NMSMAX algorithm is a frequently practiced implementation of the Nelder-Mead simplex method \cite{nelder1965simplex} and is made available in the Matrix Computation Toolbox \cite{Higham:MCT}. The NMSMAX algorithm was initialized with a regular simplex in accordance with the present algorithm.
	
	Each algorithm was given a budget of $1000$ function evaluations. Since every test function has a known global minimum, $f^*$, the tolerance of having convergence of an algorithm was $10^{-6}$, i.e. $f(\boldsymbol{x}_k) - f^* < 10^{-6}$. The internal stopping conditions of each of the algorithms were set to $0$ so that an algorithm only stops when this benchmark condition was satisfied or when it exceeded its evaluation budget. These internal stopping conditions include the gradient condition \cref{eq:tr:firstorderopt} and the minimal simplex size for NMSMAX. Only the total number of function evaluations and the difference between the final objective function value and the known minimum are recorded. We did not consider the CPU time as well as performance and data profiles which are gaining currency in the optimization community \cite{dolan2002benchmarking,more2009benchmarking}. Since the focus of the paper is primarily on the novelty of the sampling strategy and the ellipsoidal trust-region framework, we reason that the adopted approach was best suited to illustrate the elaborated methods.
	
	Each algorithm was applied using its standard settings to mimic realistic practical usage. The only parameter that was varied was the initial trust-region radius or in case of the NMSMAX algorithm, the size of the initial simplex. This corresponds with an arbitrary uniform scaling of the original optimization space which is a standard practice when performing optimization routines. The standard parameter settings of the presented algorithm are given in \cref{tab:standard settings}, those of each of the competing algorithms can be found in the respective references.
	
	We have considered $4$ different versions of the present algorithm. Details about the numerical implemention are given in the next subsection. These $4$ versions allow us to study the separate effects of the two presented mechanisms, namely the novel regular simplex based sampling strategy and the ellipsoidal trust-region mechanism. The different versions are defined in \cref{tab:variations} and are related to the values of $\mu_1 $and $\mu_2$ from \cref{algo:e2v} and \cref{eq:hessianfilter} respectively. To study the effects of the regular simplex based sampling strategy apart from the ellipsoidal trust-region framework, version V1 can be compared with the ORBIT and NMSMAX algorithm. Note that in this case the intepolation set management is purely input space oriented. In order to study the effect of the ellipsoidal trust-region framework version V2 can be compared with respect to V1. Recall that in this case the shape of the trust-region is adapted to the output behaviour of the objective function. In order to study the effect of the novel sampling strategy in combination with the ellipsoidal trust-region, versions V2 to V4 can be compared. Note that when parameter $\mu_2$ \cref{eq:hessianfilter} is equal to $1$, this corresponds with a spherical trust-region framework. It is not possible to study the behaviour of the ellipsoidal trust-region framework in combination with the sampling strategy of ORBIT, given that this strategy does not assure that a better-than-linear model is available every iteration.
	
	\begin{table}[t!]
		\caption{Standard settings of the algorithm.}
		\label{tab:standard settings}
		\centering
		\begin{tabular}{|cccccccccc|} \hline
			$\theta_1$ & $\eta_1$ & $\eta_2$ & $\gamma_1$ & $\gamma_2$ & $\beta_1$ & $\beta_2$ & $\Delta_m$ & $\Delta_M$ & $q_M$ \\ \hline\hline
			$1.25$ & $0$ & $0.6$ & $0.5$ & $2$ & $10^{-2}$ & $10^{-1}$ & $0$ & $3\Delta_0$ & $100$ \\ \hline
		\end{tabular}
	\end{table}
	
	\begin{table}[t!]
		\caption{Tested versions of the algorithm.}
		\label{tab:variations}
		\centering
		\begin{tabular}{cc|ccc} 
			$\mu_1$ &  & 0.5 & 0.75 & 1 \\
			\hline
			\multirow{2}{*}{$\mu_2$} &  0.95 & V2 & V3 & V4 \\
			& 1 & V1 & &  
			
		\end{tabular}
	\end{table}
	
	\subsection{Numerical implementation} All described procedures were implemented in the \MATLAB~ computational environment. Some of the subproblems were solved using available implemantations of the routine in \MATLAB. These are documented next.
	
	As mentioned in \cref{sec:fullylin}, once an affine subset was determined using \cref{algo:e2v}, the available set of interpolation points was reconsidered so to enhance the modelling capabilities of the surrogate model. Additional points were added as long as such did not corrupt the well-poisedness of \cref{eq:rbf-linsys} and if the total amount of interpolation points did not exceed the limit $q_M$. A similar procedure is performed by ORBIT and was adopted from that source code.
	
	The interpolation framework used here is that of Universal Kriging. We employed the DACE toolbox to obtain the modelling parameters as described in \cref{sec:RBF}. This framework requires the user to define lower and upper limits on te hyperparameter $\boldsymbol{\beta}$, respectively $\beta_1$ and $\beta_2$. More specifically we used Gaussian RBF's.
	
	Subproblem \cref{eq:TR:subproblem} was solved using the \MATLAB~ function $\mathrm{fminsearchcon}$, an extension of function $\mathrm{fminsearch}$ made available on Mathworks File exchange \cite{fminsearchcon}. Since absolute CPU time was not considered, the solution procedure was executed until fully converged.
	
	\subsection{Discussion of results}
	\subsubsection{Benchmark problems}
	By assessing the results displayed in \cref{tab:benchmark} one can observe that overall the proposed algorithm outperforms both NMSMAX and ORBIT, disregarding some coincidential draws. 
	
	Considering that V1 outperformed ORBIT in all but one occasion, we can conclude that the regular simplex based sampling strategy offers a significant advantage. One might however argue that the increased implementational and computational complexity of \cref{algo:e2v} does not justify the advantageous optimization setting in case that less challenging optimization problems are considered. In none of the benchmark problems V1 outperformed V2, V3 or V4, indicating that engaging an ellipsoidal trust-region framework is clearly a beneficial endeavour regardlessly, given the low implementational cost. The downside is however that to obtain a curved model, a minimum of $n+2$ point is required making the advanced regular simplex sampling strategy indispensable.

	Whether V2, V3 or V4 outperformed one another depends on the objective function at hand. However overall we could state that V2 is the better candidate for lower dimensional problems whilst V3 is the better candidate for higher dimensional problems as will be confirmed in the next section. Also considering the convergence history of the different algorithms it was noted that the ORBIT algorithm occosionally stagnates on a local plateau of the objective function. Since this behaviour was also exhibited by V1 (e.g. GULF, TRID; although V1 recovered each time) this demonstrates that compensating for the output behaviour of the objective function through the input space oriented sampling strategy and by altering the shape of the trust-regions is quintessential to recover from such apparent plateaus.
	
	\subsubsection{Extended Powell Singular and Extended Rosenbrock}
	
	\begin{figure}[t!]
		\centering
		\subfloat[Extended Powell Singular, $n=8$]{\label{fig:numerical:Powell}
			\includegraphics[width=0.45\textwidth]{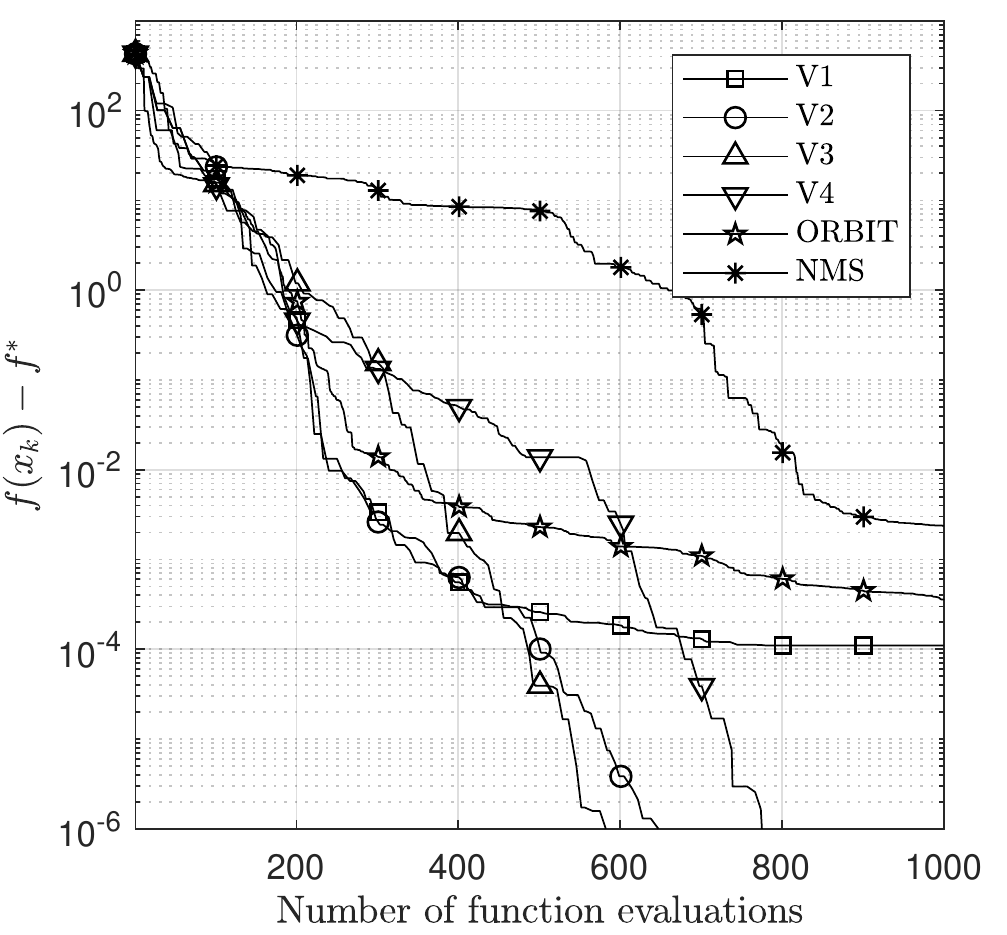}
		}
		\subfloat[Extended Rosenbrock, $n=8$]{\label{fig:numerical:Rosenbrock}
			\includegraphics[width=0.45\textwidth]{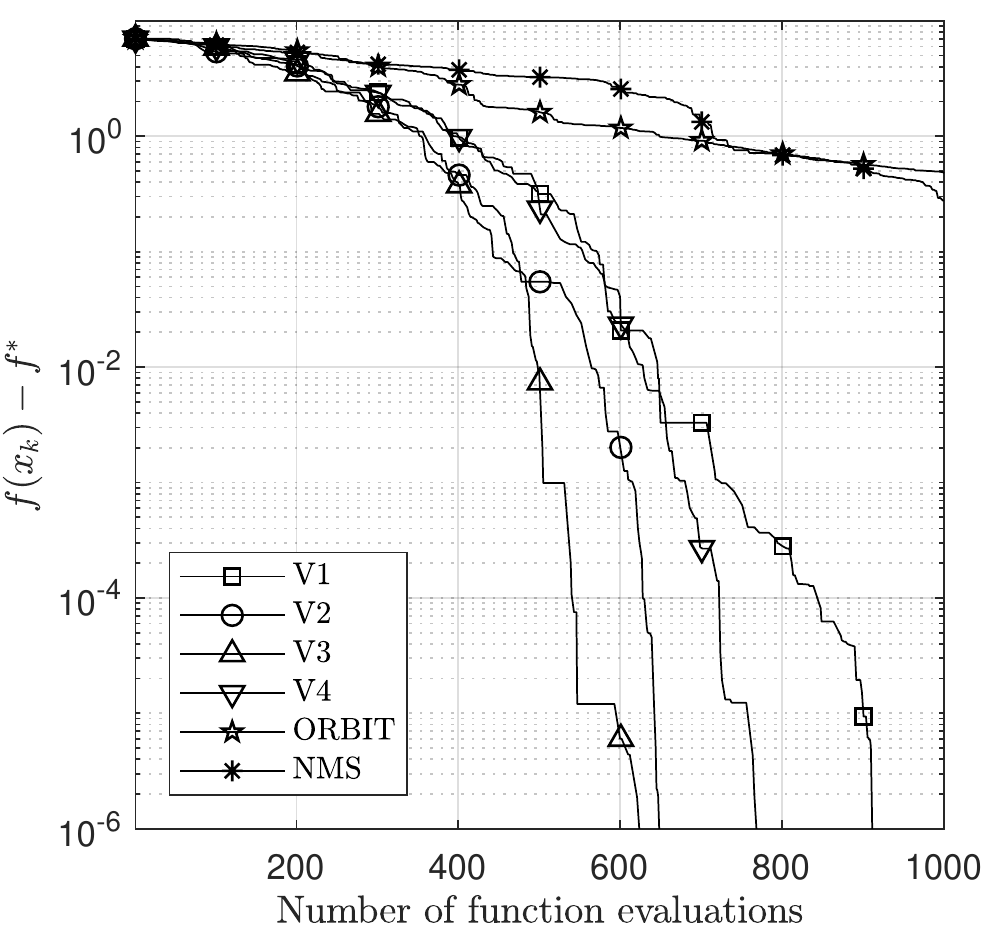}
		}
		\caption{Convergence history represented by the evolution of the best relative function value in function of the total number of function evaluations.}
		\label{fig:numerical}
		\vspace*{-.5cm}
	\end{figure}
	
	In order to demonstrate and comment on the mechanics of the elaborated methods in greater detail, convergence results on the $8$-dimensional Extended Powell Singular and Extended Rosenbrock are presented in \cref{fig:numerical} and discussed below.
	Consider the convergence history for the Extended Powell Singular function in \cref{fig:numerical:Powell}. We will initially discuss the history of the algorithms V1, ORBIT and NMSMAX. All exhibit about the same convergence behaviour during the initial stage of the optimization (up to $100$ function evaluations). After a 100 function evaluations the NMSMAX algorithm clearly starts to diverge from this trend. This is due to the fact that the NMSMAX algorithm is bounded to operate on its vertices. At the common point of about a 100 function evaluations, it reaches a region where the simplex size needs to shrink sufficiently so to attain a size that is small enough to accurately grasp the local curvature of the objective function. A similar phenomenon is observed for the ORBIT algorithm which starts to diverge from the global trend after about 200 function evaluations. The phenomenon occurs later for the ORBIT algorithm since it is not bounded to operate on the vertices of its standard simplex and the solution of \cref{eq:TR:subproblem} is only limited to the convex set, $\mathcal{B}(\boldsymbol{x}_k,\Delta_k)$. Therefore it can maintain a larger trust-region throughout the iterations than NMSMAX. That algorithm V1 can maintain a larger trust-region longer than both ORBIT and NMSMAX is a result of the improved geometry of the interpolation points. Apart from that, also algorithm V1 fails to reach the benchmark condition within the admissible budget. Now consider the behaviour of algorithm V2. The sampling strategy is equivalent to that of V1 but the modelling space is now normalized with respect to the estimated curvature of the objective function. This corresponds with an ellipsoidal trust-region geometry in the original optimization space. The effect is hardly observable during the initial stage but clearly demonstrates its usefulness after about 500 function evaluations where V2 diverges from V1 and successfully reaches the benchmark condition. It is argued that due to the transformation the modelling and optimization steps of the algorithm are now executed in a space where the objective function behaves more as if it was a quadratic function. The effect of the parameter $\mu_1$ can be assessed by comparing V2, V3 and V4. The parameter determines how strictly the volume of an exact regular simplex has to be respected. This results in a larger amount of interpolation points being added in every iteration. Therefore algorithm V3 and V4 diverge from the common trend at about the same moment that V1 and ORBIT do. In this region of the optimization space, the size of the trust-region (simplex) starts to shrink for all algorithms. This means that less function evaluations can be recycled from one trust-region to another. Algorithms V3 and V4 set a stricter limit on the volume of the simplex and are therefore forced to add more points to complement the recycled set. As a results they consume more function evaluations than V1 and V2 to cross the same region. Nonetheless this increase in function evaluations benefits the overall modelling procedure and V2 is surpassed by V3 eventually.
	
	These observations are confirmed by the convergence history for the Extended Rosenbrock function, see \cref{fig:numerical:Rosenbrock}. Here the results are directly affected by the local behaviour of the objective function and all algorithms start to diverge at the starting point. In this case V1 does reach the benchmark condition, again illustrating the improved scattering of points with respect to state-of-the-art. That as a result of the regular simplex based sampling strategy. The beneficial effect of the ellipsoidal trust-region framework is again clearly present as V2 to V4 converge faster than V1. In conclusion one can observe that the same effect is present for V3 as was the case for the Extended Powell Singular function. Function evaluations are consumed at a higher rate in the initial stage. This phenomenon ultimately benefits the modelling procedure in the final optimization stage resulting in a steeper convergence history and even the surpassing of V2 with respect to convergence speed.
	
	\begin{sidewaystable}[t!]
		\caption{Compartive results between 4 versions of the algorithm, ORBIT and NMSMAX}
		\label{tab:benchmark}
		\centering
		\begin{tabular}{lc | rrr rrr | rrr rrr} 
			Name & $n$ & \multicolumn{6}{c|}{Total number of fuction evaluations} & \multicolumn{6}{c}{Absolute difference to known minimum} \\
			& & V1 & V2 & V3 & V4 & ORBIT & NMS & V1 & V2 & V3 & V4 & ORBIT & NMS \\
			\hline\hline
			BEALE & 2 & 28 & 25 & 32 & 42 & 419 & 63 & 5,75e-07 & 6,25e-09 & 3,10e-07 & 1,34e-07 & 5,16e-10 & 2,93e-07 \\ 
			BRANIN & 2 & 24 & 25 & 32 & 36 & 53 & 56 & 1,10e-07 & 5,35e-07 & 7,60e-08 & 6,14e-07 & 1,10e-05 & 7,81e-07 \\ 
			ROSENBROCK & 2 & 37 & 37 & 33 & 39 & 838 & 100 & 2,51e-07 & 9,56e-08 & 9,40e-07 & 1,69e-07 & 1,80e-03 & 3,71e-07 \\ 
			SCHWEFEL36 & 2 & 26 & 25 & 26 & 32 & 69 & 69 & 9,35e-07 & 4,29e-07 & 1,44e-07 & 7,43e-08 & 9,17e-07 & 4,39e-07 \\ 
			GULF & 3 & 494 & 159 & 165 & 226 & 1000 & 389 & 3,29e-07 & 7,71e-07 & 1,84e-09 & 8,95e-07 & 2,66e+00 & 5,60e-07 \\ 
			HART3 & 3 & 52 & 46 & 49 & 69 & 348 & 129 & 9,98e-09 & 8,33e-07 & 7,14e-07 & 2,64e-08 & 3,91e-06 & 6,08e-07 \\ 
			HELIX & 3 & 53 & 42 & 48 & 48 & 443 & 115 & 8,29e-07 & 5,97e-07 & 3,08e-07 & 4,50e-07 & 4,61e-04 & 4,31e-07 \\ 
			ROSENBROCK & 3 & 77 & 68 & 80 & 88 & 1000 & 176 & 3,81e-07 & 4,41e-07 & 6,66e-07 & 1,18e-07 & 5,20e-03 & 5,12e-07 \\ 
			SCHWEFEL4 & 3 & 55 & 46 & 41 & 67 & 77 & 93 & 1,72e-07 & 7,98e-08 & 7,20e-07 & 1,81e-08 & 8,05e-07 & 7,07e-07 \\ 
			BROWNDEN & 4 & 118 & 84 & 112 & 121 & 398 & 236 & 6,77e-07 & 9,34e-06 & 5,16e-06 & 3,44e-06 & 9,60e-06 & 9,77e-06 \\ 
			HART4 & 4 & 50 & 42 & 48 & 38 & 62 & 122 & 3,86e-07 & 2,78e-07 & 3,95e-07 & 2,94e-07 & 5,56e-07 & 8,85e-07 \\ 
			POWELL & 4 & 201 & 92 & 104 & 148 & 615 & 366 & 6,77e-07 & 7,43e-07 & 8,85e-07 & 4,51e-07 & 6,01e-05 & 9,45e-07 \\ 
			ROSENBROCK & 4 & 136 & 111 & 143 & 144 & 1000 & 332 & 8,67e-07 & 9,90e-07 & 3,10e-07 & 7,74e-07 & 2,75e-02 & 8,78e-07 \\ 
			SCHWEFEL4 & 4 & 63 & 61 & 62 & 74 & 112 & 156 & 2,46e-07 & 9,04e-07 & 5,55e-07 & 3,18e-07 & 9,25e-07 & 9,75e-07 \\ 
			WOOD & 4 & 166 & 79 & 85 & 109 & 255 & 217 & 6,35e-07 & 5,16e-07 & 3,78e-07 & 2,13e-08 & 4,04e-04 & 7,57e-07 \\ 
			BIGGS & 6 & 660 & 294 & 243 & 529 & 156 & 1000 & 7,54e-04 & 1,52e-07 & 9,65e-07 & 9,48e-07 & 6,60e-03 & 5,68e-05 \\ 
			HART6 & 6 & 104 & 116 & 95 & 107 & 139 & 179 & 5,63e-07 & 5,43e-07 & 5,23e-07 & 4,08e-07 & 8,21e-07 & 9,39e-07 \\ 
			ROSENBROCK & 6 & 302 & 257 & 264 & 435 & 1000 & 818 & 3,02e-07 & 6,78e-07 & 8,09e-07 & 7,34e-07 & 7,22e-02 & 9,90e-07 \\ 
			SCHWEFEL4 & 6 & 124 & 116 & 111 & 169 & 173 & 221 & 8,37e-07 & 7,10e-07 & 2,38e-07 & 7,37e-07 & 9,47e-07 & 8,81e-07 \\ 
			TRID & 6 & 405 & 206 & 205 & 194 & 399 & 541 & 6,30e-03 & 6,67e-07 & 2,77e-07 & 1,79e-07 & 6,90e-03 & 7,94e-07 \\ 
			WATSON & 6 & 117 & 83 & 97 & 119 & 222 & 305 & 2,11e-07 & 8,63e-07 & 1,56e-07 & 4,16e-07 & 9,83e-07 & 7,21e-07 \\ 
			POWELL & 8 & 168 & 135 & 139 & 151 & 225 & 410 & 7,61e-07 & 2,88e-07 & 5,77e-07 & 3,95e-07 & 1,00e-06 & 9,70e-07 \\ 
			ROSENBROCK & 8 & 912 & 648 & 624 & 769 & 1000 & 1000 & 7,41e-07 & 7,25e-07 & 9,35e-07 & 8,96e-07 & 4,81e-01 & 2,75e-01 \\ 
			TRID & 8 & 1000 & 648 & 585 & 778 & 1000 & 1000 & 1,10e-04 & 9,81e-07 & 8,49e-07 & 3,68e-07 & 3,53e-04 & 2,40e-03 \\ 
		\end{tabular}
	\end{sidewaystable}
	
	\section{Conclusion}
	\label{sec:conclusions}
	This paper is on reducing the number of objective function evaluations in derivative-free numerical optimization with surrogate models. Radial basis functions are here considered as surrogate models with the interpolation framework provided by Universal Kriging. In each iteration of an interpolation based optimization algorithm, a trust region can be used wherein an optimization is performed using the surrogate model. The interpolation framework consists of using an interpolation set of which a subset of interpolation points needs to be well-poised so that the quality of the surrogate model is maximal. Interpolation set management aims at maximizing this range of model validity.
	
	Currently, at least $q=n+1$ interpolation points are considered to guarantee a valid surrogate model. We advanced on interpolation based optimization by relying on a subset of at least $q=n+2$ interpolation points within the interpolation set management to improve the scattering of the interpolation points but more importantly to have surrogate models that exhibit curvature. For that purpose we defined fully linear surrogate models in \cref{def:fullylinear} that can be incorporated in the trust-region algorithm (\cref{alg:TR:fullylinear}). We implemented RBF models using a flexible multivariate interpolation framework for which Universal Kriging is used (\cref{sec:RBF}). In \cref{lemma:bthanlinear} we provide for a given RBF model with linear tail, information on the boundedness of the tolerance of the RBF model and objective function as well on their gradients, i.e. constants $\kappa_f$ and $\kappa_g$ in \cref{def:fullylinear}. This way we are able to predict whether an interpolation model classifies as being fully linear. Since we consider $q=n+2$ interpolation points that represent a regular simplex \cref{sec:sub:regular-simplex}, we propose in \cref{algo:e2v} and \cref{algo:se} how the interpolation set management needs to be adapted, i.e. how the input points that need to be evaluated in the objective function need to be chosen.
	
	By leveraging on the curvature of the surrogate model that is valid in the presented interpolation set management we introduced an ellipsoidal trust region framework that allows morphed trust-regions that adapt to the objective function geometry. This adaptation is implemented by means of curvature normalization with coordinate transformation \cref{eq:affinetrans} and filtered Hessian approximation \cref{eq:hessianfilter}. This way we allow to incorporate the output behavior of the objective function in the shape of the trust-region.
	
	The proposed interpolation set management with proper input of interpolation points that follow a well-poised geometry together with the output behavior of the objective function that is incorporated in the ellipsoidal trust region, have shown their benefits with respect to the reduction of the number of objective function evaluations when performing elaborate numerical optimizations on various test functions with $n\leq8$. In order to show the impact of internal parameters of the algorithms, i.e. $\mu_1$ in \cref{algo:e2v} that is related to the input set and $\mu_2$ in \cref{eq:hessianfilter} that is related to the ellipsoidal trust region shape, results of 4 different versions of the presented algorithm are shown in \cref{fig:numerical} and \cref{tab:benchmark}. Results clearly illustrate the impact of the used input interpolation points and the ellipsoidal trust region. Moreover, the convergence histories show faster convergence compared to algorithms incorporating at least $n+1$ interpolation points and that use a spherical trust region.

	\section*{Acknowledgements}
	The authors acknowledge the support of the FWO research project G.0D93.16N and the Flanders Make project EMODO.
	
	\newpage
	\appendix
	\section{Universal Kriging}
	\label{appendix:unikrig}
	The Universal Kriging (UK) framework obtains an interpolating approximation model for the function $f:\mathcal{X}\subset\mathbb{R}^n\rightarrow f(\mathcal{X})\subset\mathbb{R}$ as the {most likely} unbiased linear predictor, $m(\boldsymbol{x})$, based on a set of $q$ interpolation points $\mathcal{Z}=\{\boldsymbol{z}_1,\dots,\boldsymbol{z}_q\}, \boldsymbol{z}_i\in\mathcal{X}$, and, a set of $q$ function evaluations, $\mathcal{F} = \{f(\boldsymbol{z}_1),\dots,f(\boldsymbol{z}_q)\}$. The framework presumes that the modelled function $f$ can be decomposed as 
	\begin{equation}
	f(\boldsymbol{x}) = p(\boldsymbol{x}) + w(\boldsymbol{x}) = \boldsymbol{\pi}(\boldsymbol{x})^\top\boldsymbol{\beta} + w(\boldsymbol{x})
	\end{equation}
	
	The function $p(\boldsymbol{x})$ is assumed to be an element of a linear function space, $\mathcal{A}$. In the UK framework one is interested in the specific case where $\mathcal{A}$ is the $n$-variate polynomial space, $\mathcal{P}^n_d$, of at most degree $d$ spanned by the arbitrary basis, $\hat{\mathcal{P}}^n_{d} = \{\pi_1(\boldsymbol{x}),\dots,\pi_{\hat{d}}(\boldsymbol{x})\}$. We further introduce the vector, $\boldsymbol{\pi}(\boldsymbol{x})$, which entries are determined by $\pi_i(\boldsymbol{x})$, so that we can write, $p(\boldsymbol{x}) = \pi(\boldsymbol{x})^\top\boldsymbol{\beta}$, where $\boldsymbol{\beta}\in\mathbb{R}^{\hat{d}}$ is a unique coefficient vector. It is further assumed that $w(\boldsymbol{x})$ is an unbiased random process and that the output values $w(\boldsymbol{x})$ and $w(\boldsymbol{z})$ are correlated. Their correlation is modelled by a function from the class of radial basis functions, i.e. $\mathrm{corr}(w(\boldsymbol{x}),w(\boldsymbol{z})) = \phi(\|\boldsymbol{x}-\boldsymbol{z}\|)$. One looks for an unbiased linear predictor, $m(\boldsymbol{x}) = \boldsymbol{\eta}(\boldsymbol{x})^\top \boldsymbol{f}$, where he elements of $\mathcal{F}$ determine the entries of $\boldsymbol{f}$.
	
	The model decomposition can be applied to this vector $\boldsymbol{f}$ so that one obtains $\boldsymbol{f} = {\boldsymbol{\Pi}}\boldsymbol{\beta} + \boldsymbol{w}$. Here the entries of $\boldsymbol{w}$ are defined as $w(\boldsymbol{x}_i) = f(\boldsymbol{x}_i)-p(\boldsymbol{x}_i)$ and the matrix ${\boldsymbol{\Pi}}\in\mathbb{R}^{q\times\hat{d}}$ as $\Pi_{ij} = \pi_j(\boldsymbol{x}_i)$. Substitution in the linear predictor yields that $m(\boldsymbol{x}) = \boldsymbol{\eta}(\boldsymbol{x})^\top {\boldsymbol{\Pi}} \boldsymbol{\beta} + \boldsymbol{\eta}(\boldsymbol{x})^\top \boldsymbol{w}$. Hence, we can define an error function between the linear predictor and the true function as ${\Delta}(\boldsymbol{x}) = m(\boldsymbol{x}) - f(\boldsymbol{x}) = (\boldsymbol{\eta}(\boldsymbol{x})^\top\boldsymbol{\Pi} - \boldsymbol{\pi}(\boldsymbol{x})^\top)\boldsymbol{\beta} + (\boldsymbol{\eta}(\boldsymbol{x})^\top \boldsymbol{w} - w(\boldsymbol{x}))$.
	
	We can now impose two additional properties on the linear predictor that affect the error function $\Delta$ and that will render the predictor unique. Firstly, the linear predictor should be unbiased, which can be expressed as $\mathrm{E}[\Delta] \equiv 0$. It follows that ${\boldsymbol{\Pi}}^\top\boldsymbol{\eta} = \boldsymbol{\pi}$ and consequently that $\Delta = \boldsymbol{\eta}^\top \boldsymbol{w} - w$. Secondly, we wish to determine $\boldsymbol{\eta}$ so that the expected squared error is minimized. Given that $\mathrm{E}[\Delta] = 0$, the expected squared error can be defined as $\mathrm{E}[\Delta^2]$ and after substitution of the error function we obtain $\mathrm{E}[(w-\boldsymbol{\eta}^\top \boldsymbol{w})^2]$. Evaluation of the expectation operator determines that the error is proportional to the functional $\epsilon[\boldsymbol{\eta}] = \boldsymbol{\eta}^\top{\boldsymbol{\Phi}}\boldsymbol{\eta} - 2\boldsymbol{\eta}^\top\boldsymbol{\phi} + 1$, where matrix ${\boldsymbol{\Phi}}\in\mathbb{R}^{q\times q}$ is defined as $\Phi_{ij} = \phi(\|\boldsymbol{z}_i-\boldsymbol{z}_j\|)$ and the entries of $\boldsymbol{\phi}(\boldsymbol{x})$ are determined by $\phi(\|\boldsymbol{x}-\boldsymbol{z}_i\|)$. The most likely unbiased linear predictor $\boldsymbol{\eta}$ is thus the solution of the equality constrained variational problem
	\begin{equation}
	\min_{\boldsymbol{\eta},\boldsymbol{\xi}} ~ \tfrac{1}{2}\left(\boldsymbol{\eta}^\top{\boldsymbol{\Phi}}\boldsymbol{\eta} - 2\boldsymbol{\eta}^\top\boldsymbol{\phi} + 1\right) + \boldsymbol{\xi}^\top\left(\boldsymbol{\Pi}^\top\boldsymbol{\eta} - \boldsymbol{\pi}\right)
	\end{equation}
	
	Expression of the first order optimality condition yields the system of linear equations given in \cref{eq:kriging-linsys}.
	
	According to the probabilistic interpretation of the random process $w(\boldsymbol{x})$, the function values in $\boldsymbol{f}$ can be associated to the outcome of a stochastic experiment. Consequently it is possible to associate a probability value to the outcome of this experiment. If we make the additional assumption that the random process $w(\boldsymbol{x})$ itself behaves as a Gaussian with expectation $\boldsymbol{\mu}={\boldsymbol{\Pi}}\boldsymbol{\beta}$ and covariance $\sigma$, the probability of the experiment is given by
	\begin{equation}
	P(\sigma^2) = (2\pi\sigma^2)^{-\frac{n}{2}}\det({\boldsymbol{\Phi}})^{-\frac{1}{2}} \exp\left(-\tfrac{1}{2}\sigma^{-2}(\boldsymbol{f}-\boldsymbol{\mu})^\top{\boldsymbol{\Phi}}^{-1}(\boldsymbol{f}-\boldsymbol{\mu})\right)
	\end{equation}
	
	A maximum log likelihood estimation for the parameter $\sigma^2$ results in the following estimate $\hat{\sigma}^2 = \tfrac{1}{n}({\boldsymbol{f}}-{\boldsymbol{\Pi}}\boldsymbol{\beta})^\top{\boldsymbol{\Phi}}^{-1}({\boldsymbol{f}}-{\boldsymbol{\Pi}}\boldsymbol{\beta})$. Now, assuming that the correlation model $\mathrm{corr}(\boldsymbol{x},\boldsymbol{z})$ still depends on the additional hyperparameter, say $\boldsymbol{\theta}$ (e.g. $\phi(\|\boldsymbol{x}-\boldsymbol{z}\|;\boldsymbol{\theta})$). The optimal parameter $\boldsymbol{\theta}^*$ can then also be determined by maximizing the likelihood \begin{equation*}
	\boldsymbol{\theta}^* = \arg\min_{\theta} P(\hat{\sigma}^2)
	\end{equation*}
	\newpage
	\section{Error bounds on better-than-linear RBF models}
	\label{appendix:btl}
	Using the Einstein summation convention, we can write an RBF model with linear polynomial tail as presented in \cref{eq:app2:rbf}. Here we have implicitly defined the error functions $e_f(\boldsymbol{x})$ and $\boldsymbol{e}_{\nabla f}(\boldsymbol{x})$. We are interested in bounds on the values these functions can take given that the interpolation set, $\mathcal{Z}$, contains the origin and further that both $\boldsymbol{x}$ and $\mathcal{Z}$ are contained within the ball $\mathcal{B}(\boldsymbol{0},\Delta)$ 
	\begin{equation}
	\label{eq:app2:rbf}
	\begin{aligned}
	m(\boldsymbol{x}) &\doteq \gamma+\boldsymbol{x}^\top\boldsymbol{\beta}+\alpha_j\phi_j(\boldsymbol{x})= f(\boldsymbol{x}) + e_f(\boldsymbol{x}) \\
	\nabla m(\boldsymbol{x}) &\doteq \boldsymbol{\beta}+\alpha_j\nabla\phi_j(\boldsymbol{x})= \nabla f(\boldsymbol{x}) + \boldsymbol{e}_{\nabla f}(\boldsymbol{x})
	\end{aligned}
	\end{equation}
	
	Substracting the left- and right-hand side expressions for $m(\boldsymbol{x})$ from the left- and right-hand side expressions of each of the $q> n+1$ interpolation conditions, $\gamma + \boldsymbol{x}^\top\boldsymbol{\beta} + \alpha_j\phi_j(\boldsymbol{x}) = f(\boldsymbol{x}_i)$, one obtains the following set of equalities
	\begin{equation}
	\left\lbrace\begin{aligned}
	-\boldsymbol{x}^\top \boldsymbol{\beta} &= f(\boldsymbol{0}) - f(\boldsymbol{x}) - \alpha_j \left(\phi_j(\boldsymbol{0})-\phi_j(\boldsymbol{x})\right) - e_f(\boldsymbol{x}) \\
	(\boldsymbol{z}_i-\boldsymbol{x})^\top \boldsymbol{\beta}  &= f(\boldsymbol{z}_i) - f(\boldsymbol{x}) - \alpha_j \left(\phi_j(\boldsymbol{z}_i)-\phi_j(\boldsymbol{x})\right) - e_f(\boldsymbol{x}), ~\forall j > 0
	\end{aligned}\right.
	\label{eq:app:equal}
	\end{equation}
	
	These inequalities can be manipulated so to relate the error on the function value, $e_f(\boldsymbol{x})$, to the error on the derivative value, $\boldsymbol{e}_{\nabla f}(\boldsymbol{x})$. For that purpose, we shall first rewrite the functions $f(\boldsymbol{x})$ and $\phi_k(\boldsymbol{x})$ using the integral representation $f(\boldsymbol{x}) = f(\boldsymbol{z}_i)+\int_{\boldsymbol{z}_i}^{\boldsymbol{x}}\nabla f(\boldsymbol{s})^\top\text{d}\boldsymbol{s}$. After substituting $\boldsymbol{x}(1-t)+\boldsymbol{z}_i t$ for $\boldsymbol{s}$, it follows that $\int_{\boldsymbol{z}_i}^{\boldsymbol{x}}\nabla f(\boldsymbol{s})^\top\text{d}\boldsymbol{s} = -\int_{0}^{1}\nabla f(\boldsymbol{x}-t(\boldsymbol{z}_i-\boldsymbol{x}))^\top(\boldsymbol{z}_i-\boldsymbol{x})\text{d}t$. Now in order to make the derivative error function emerge, one can add $-\alpha_j\boldsymbol{x}^\top\nabla \phi_j(\boldsymbol{x})+\boldsymbol{x}^\top\nabla f(\boldsymbol{x})$ to each side of the expressions. Together with the left-hand side of \cref{eq:app:equal} this expression combines into $\boldsymbol{e}_{\nabla f}(\boldsymbol{x})$. In the right-hand side these expressions can be incorporated in the integrandum since their value is constant with respect to the integration variable and considering that the integral is taken over the unit interval. 
	
	These manipulations allow to reorganize the equalities in \cref{eq:app:equal} and rewrite them in function of the operators $\mathcal{Q}[f](\boldsymbol{x})$ and $\mathcal{Q}_i[f](\boldsymbol{x})$, which are defined as
	\vspace*{-5pt}
	\begin{equation*}
	\left\lbrace\begin{aligned}
	\mathcal{Q}[f](\boldsymbol{x}) &= \int_{0}^{1}\left(\nabla f(\boldsymbol{x}-t\boldsymbol{x})-\nabla f(\boldsymbol{x})\right)^\top(-\boldsymbol{x})\text{d}t \\ \mathcal{Q}_i[f](\boldsymbol{x}) &= \int_{0}^{1}\left(\nabla f(\boldsymbol{x}+t(\boldsymbol{z}_i-\boldsymbol{x}))-\nabla f(\boldsymbol{x})\right)^\top(\boldsymbol{z}_i-\boldsymbol{x})\text{d}t
	\end{aligned}\right.
	\vspace*{-5pt}
	\end{equation*}
	and in function of the error functions, $e_f(\boldsymbol{x})$, and, $\boldsymbol{e}_{\nabla f}(\boldsymbol{x})$, as given by
	\vspace*{-5pt}
	\begin{equation}
	\left\lbrace\begin{aligned}
	-\boldsymbol{x}^\top \boldsymbol{e}_{\nabla f}(\boldsymbol{x}) &= \mathcal{Q}[f](\boldsymbol{x}) - \alpha_k \mathcal{Q}[\phi_k](\boldsymbol{x}) - e_f(\boldsymbol{x}) \\
	(\boldsymbol{z}_i-\boldsymbol{x})^\top \boldsymbol{e}_{\nabla f}(\boldsymbol{x}) &= \mathcal{Q}_i[f](\boldsymbol{x}) - \alpha_j \mathcal{Q}_i[\phi_j](\boldsymbol{x}) - e_f(\boldsymbol{x}), ~\forall i > 0
	\end{aligned}\right.
	\label{eq:derivativelinkvalue}
	\end{equation}
	
	Introducing the operator, $\delta\mathcal{Q}_i[f](\boldsymbol{x})= \mathcal{Q}_i[f](\boldsymbol{x}) - \mathcal{Q}[f](\boldsymbol{x})$, and substracting the first equation from the others yields the final expressions
	\begin{equation}
	\label{eq:errortooperator}
	\boldsymbol{z}_i^\top \boldsymbol{e}_{\nabla f}(\boldsymbol{x}) = \delta \mathcal{Q}_i[f](\boldsymbol{x}) - \alpha_k \delta \mathcal{Q}_i[\phi_k](\boldsymbol{x}),  ~\forall i > 0
	\end{equation}
	
	Considering that $\mathcal{Z}\subset\mathcal{B}(\boldsymbol{0},\Delta)$, it can be shown that the right hand side is bounded by (see \cite{Conn2008,dennis1996numerical})
	\begin{equation}
	\left\lbrace\begin{aligned}
	\left|\mathcal{Q}[f](\boldsymbol{x}) - \alpha_j \mathcal{Q}[\phi_j](\boldsymbol{x})\right| &\leq \tfrac{1}{2}(\gamma_f+\gamma_\phi)\Delta^2 \\
	\left|\mathcal{Q}_i[f](\boldsymbol{x}) - \alpha_j \mathcal{Q}_i[\phi_j](\boldsymbol{x})\right| &\leq 2(\gamma_f+\gamma_\phi)\Delta^2 \\
	\left|\delta \mathcal{Q}_i[f](\boldsymbol{x}) - \alpha_j \delta\mathcal{Q}_i[\phi_j](\boldsymbol{x})\right| &\leq \tfrac{5}{2}(\gamma_f+\gamma_\phi)\Delta^2 
	\end{aligned}\right.
	\end{equation}
	where $\gamma_f$ and $\gamma_\phi$ are Lipshitz constants for $\nabla f(\boldsymbol{x})$ and $\alpha_k\nabla\phi_k(\boldsymbol{x})$ given that $\boldsymbol{x}\in\mathcal{B}(\boldsymbol{0},\Delta)$.
	
	Combination of the equations in \cref{eq:errortooperator} yields a bound on the derivative error function in function of the interpolation points contained in $\mathcal{Z}$
	\begin{equation*}
	\left\|\begin{matrix}
	\boldsymbol{z}_1^\top \\
	\vdots\\
	\boldsymbol{z}_q^\top
	\end{matrix}\right\| \|\boldsymbol{e}_{\nabla f}(\boldsymbol{x})\| = \left\|\mathbf{Z}\right\|  \|\boldsymbol{e}_{\nabla f}(\boldsymbol{x})\| \ \leq \tfrac{5}{2}\sqrt{q-1} (\gamma_f+\gamma_\phi) \Delta^2
	\end{equation*}
	
	Since we already had established a relation between both error functions in \cref{eq:derivativelinkvalue}, the bound above allows one to determine a bound on $e_f(\boldsymbol{x})$ as well, that is given by
	\begin{equation}
	|e_f(\boldsymbol{x})| \leq \frac{1}{2}(\gamma_f+\gamma_\phi)\Delta^2 + \tfrac{5}{2}\sqrt{q-1}\left\|\mathbf{Z}\right\|^{-1} (\gamma_f+\gamma_\phi) \Delta^3
	\end{equation}
	
	Now introducing the scaled interpolation matrix $\Delta\|\hat{\mathbf{Z}}\|=\left\|\mathbf{Z}\right\|$, we have established taylor-like bounds on both the error value functions that are related to the inverse value of $\|\hat{\mathbf{Z}}\|$. These bounds are in fact solely determined by the interpolation set geometry.
	\begin{equation}
	\left\lbrace\begin{aligned}
	\|\boldsymbol{e}_{\nabla f}(\boldsymbol{x})\| &\leq \tfrac{5}{2}\sqrt{q-1} \|\hat{\mathbf{Z}}\|^{-1} (\gamma_f+\gamma_\phi) \Delta \\
	|\boldsymbol{e}_{f}(\boldsymbol{x})| &\leq \left(\tfrac{5}{2}\sqrt{q-1} \|\hat{\mathbf{Z}}\|^{-1} + \tfrac{1}{2} \right) (\gamma_f+\gamma_\phi) \Delta^2
	\end{aligned}\right.
	\end{equation}
	
	\newpage
	\section{Some simplex related defenitions}
	\label{appendix:simplex}
	Consider a set of $n+1$ points $\mathcal{Z}=\{\boldsymbol{z}_0,...,\boldsymbol{z}_n\}$ associated to a set of coordinates in an $n$ dimensional space. Assume that these points are affinely independent so that also the column vectors $\{\boldsymbol{z}_1-\boldsymbol{z}_0,...,\boldsymbol{z}_n-\boldsymbol{z}_0\}$ are linearly independent. A simplex is therefore the $n$-dimensional generalization of a triangle. In simplex terminology one has that each point from $\mathcal{Z}$ determines a vertice, each couple of vertices determines an edge and each subset of $n$ vertices determines a face. A $n$-simplex thus has $n+1$ vertices, $\tfrac{n(n+1)}{2}$ edges and $n+1$ faces. Conveniently it follows that each face of an $n+1$-simplex constitutes an $n$-simplex itself.
	
	According to the general convention, when speaking of the $n$-simplex that is determined by $\mathcal{Z}$ one reffers to the convex set determined by $\mathcal{S}$
	\begin{equation}
	\mathcal{S} = \left\{\sum_{i=0}^n \theta_i\boldsymbol{z}_k\left|\sum_{i=0}^n\theta_i=1\right. \wedge\theta_i \geq 0,\forall i\right\}
	\end{equation}
	
	Now as it may be inconsistent with this general convention, in this work we refer to the set of vertices $\mathcal{Z}$ when speaking of the simplex. However, when we speak of the hypervolume of that simplex, we do in fact refer to the hypervolume of the convex set, $\mathcal{S}$.
	
	The volume $V$ of the simplex determined by $\mathcal{Z}$ is defined as the determinant of the associated column vector matrix $[\mathcal{Z}]$ normalized by the factor $n!$
	\begin{equation}
	\label{eq:simvolume}
	V(\mathcal{Z}) = \frac{1}{n!}\det([\mathcal{Z}])
	\end{equation}
	here the operator $[\cdot]$ is defined as $
	[\mathcal{Z}]=\begin{bmatrix}
	\boldsymbol{z}_1-\boldsymbol{z}_0&\cdots&\boldsymbol{z}_n-\boldsymbol{z}_0
	\end{bmatrix}
	$.
	
	In the following subsection we mention some relevant simplex geometries and properties.
	
	\subsection{Orthocentric simplex}
	A simplex for which the following holds are called orthocentric simplices: {each edge is perpendicular to all the edges it does not meet} \cite{thealtitudes}. Orthocentric simplices possess a unique orthocentre. In the special cases where $n\leq2$, we have that every set of points that determine a simplex are orthocentric as well. 
	
	In \cite{gerber1975}, it is shown that an orthocentric simplex maximizes the volume satisfying a predetermined orthocentre, $\boldsymbol{o}$, and, distances from the simplex vertices to the orthocentre $d_i=\|\boldsymbol{z}_i-\boldsymbol{o}\|$. These conditions determine the simplex, $\mathcal{Z}$, within an isometry. It is this very property that we exploit to maximize the volume of the union of an existing set and an expansion set generated by the procedure $\mathrm{Simexpand}$, \cref{algo:se}. The expansion set is constructed exactly so that all of its edges are perpendicular to the edges determined by the existing set. Given the freedom of the edge lengths, beit that the union of exesting and expansion set is restricted to an $n$-sphere, these are determined so to maximize the volume.
	
	\subsection{Standard simplex}\label{appendix:simplex:standard} A standard $n$-simplex is obtained when $n$ of its edges are orthogonal and have the same length, $l$. A convenient coordinate representation of a standard $n$-simplex is given by the origin, $\boldsymbol{0}$, and the scaled coordinates of the $n$ basis vectors, $l\boldsymbol{e}_i$. It is easily verified that the $\tfrac{n(n-1)}{2}$ non-orthogonal edges have length $\sqrt{2}l$. 
	
	We also note that the standard simplex is not an orthocentric simplex.
	
	\subsection{Regular simplex}
	A regular $n$-simplex with side length $l$ is obtained when all edges have the same length. From this definition it follows immediately that any subset of $m+1$ points from that simplex determines a regular $m$-simplex.
	
	If we consider all simplices $\mathcal{Z}$ which vertices are contained within an $n$-hypersphere with radius $\Delta$, it can be shown that if a simplex's volume is maximized it is an element from the class of regular simplices \cite{thealtitudes}. Moreover it can be shown that the regular simplex is an orthocentric simplex.
	
	It is worth mentioning that a regular $n$-simplex is embedded in a standard $n+1$-simplex as the face determined by the $\tfrac{n(n+1)}{2}$ non-orthogonal edges of the standard $n+1$-simplex. The edge length of a regular $n$-simplex that is embedded in a standard $n+1$-simplex with orthogonal edge length, $l$, is $\sqrt{2}l$ in accordance with the non-orthogonal edge length of a standard simplex with edge lenght $l$ as determined above.
	
	\subsubsection{Volume of the regular $n$-simplex}
	\label{apendix:volume}
	Let $V_{n}^s(l)$ be the hypervolume occupied by the standard $n$-simplex with orthogonal edge length $l$. Following the definition given in \cref{eq:simvolume}, it is easily verified that $V_{n+1}^s(l) = \tfrac{1}{(n+1)!}l^{n+1}$. Now consider that the face determined by all the non-orthogonal edges of the standard $n+1$-simplex constitutes a regular $n$-simplex. The distance from the intersection of the orthogonal edges, i.e. origin $\boldsymbol{0}$ according the coordinate representation in \ref{appendix:simplex:standard}, to the face that determines the regular $n$-simplex is given by $\tfrac{\sqrt{n+1}}{n+1}l$. Since the hypervolume of a standard $n$-simplex can also be determined using the definition of the hyperpyramid, i.e. by multiplying the hypervolume of one of its faces with the corresponding altitude divided by $n$, the hypervolume of the regular $n$-simplex $V_n^r(l)$ can be determined. Recall that the edge length of the regular $n$-simplex determined by the non-orthogonal edges is $\sqrt{2}l$ and that therefore the result must be scaled by a factor $\sqrt{2^n}$.
	\begin{equation}
	\frac{1}{n+1}\cdot\frac{\sqrt{n+1}}{n+1}l\cdot V_n^r(\sqrt{2}l) = V_{n+1}^s(l)\Rightarrow  V_n^r(l)=\frac{\sqrt{n+1}}{n!}\frac{l^n}{\sqrt{2^n}}
	\end{equation}
	
	%
	
	\subsubsection{Radius of the circumscribed $n$-sphere} The radius, $r_n(l)$, of the $n$-sphere circumscribing a regular $n$-simplex with side length $l$ can be determined reconsidering the geometric relations between a standard $n+1$-simplex and a regular $n$-simplex described to determine the hypervolume.  In a regular $n$-simplex the distance between the unique orthocentre and a vertice is equal to the radius of the circumscribed $n$-sphere. It holds that the altitude of the face determined by the non-orthogonal edges and one of the orthogonal edges of the standard $n+1$-simplex constitute a right triangle with the orthocentre of the associated regular $n$-simplex and a common vertice.  Note that again the edge length of this $n$-simplex is $\sqrt{2}l$ and that therefore the result must be scaled with by a factor $\sqrt{2}$.
	\begin{equation}
	l^2 = r_n(\sqrt{2}l)^2 + \frac{1}{n+1}l^2 \Rightarrow r_n(l) = \frac{\sqrt{n}}{\sqrt{n+1}}\frac{l}{\sqrt{2}}
	\end{equation}
	
	The side length, $l_n(\Delta)$, of the regular $n$-simplex inscribed in a $n$-sphere with radius $\Delta$ is hence given by $l_n(\Delta) = \tfrac{\sqrt{2}\sqrt{n+1}}{\sqrt{n}}\Delta$
	
	\subsubsection{Other usefull relations}
	From these results we can derive some other equalities that are of interest for the procedure $\mathrm{Simexpand}$ described by \cref{algo:se}.
	
	The radius of the circumscribed $m$-sphere of a regular $m$-simplex that is part of a regular $n$-simplex circumscribed by an $n$-sphere with radius $\Delta$ is given by
	\begin{equation}
	r_m^n(\Delta)= r_m(l_n(\Delta)) = \sqrt{\tfrac{m(n+1)}{(m+1)n}}\Delta
	\end{equation}
	
	The distance, $d^n_m(\Delta)$, from the orthocentre of the regular $n$-simplex to a subset of $m+1$ vertices determining such a regular $m$-simplex, can hence be determined considering that
	\begin{equation}
	\Delta^2 = d^n_m(\Delta)^2 + r^n_m(\Delta)^2 \Rightarrow d^n_m(\Delta) = \sqrt{\tfrac{n-m}{(m+1)n}}\Delta
	\end{equation}
	
	This relation can be generalised. Consider therefore a regular $m$-simplex embedded in a $n$-dimensional space and inscribed in an $n$-sphere with radius $\Delta$. Its edge length must then be $l_m(\Delta)$. Now if one would translates this $m$-simplex along any direction perpundicalar to the $m$-dimensional subspace spanned by its edges over a distance $\sigma \Delta$, this $m$-simplex remains on the surface of the $n$-sphere as long as it is scaled by a factor $\sqrt{1-\sigma^2}\Delta$.

	%

%
%
%
	
	\newpage
	
	\bibliographystyle{unsrt}
	\bibliography{references}
	
\end{document}